\documentclass[11pt]{article}
\usepackage[T1]{fontenc}
\usepackage[margin=1in]{geometry}
\usepackage{graphicx} 
\usepackage{amsmath}
\usepackage{amssymb}
\usepackage{mathtools}
\usepackage{mathrsfs}
\usepackage{xcolor}
\usepackage{hyperref,cleveref}
\usepackage{amsthm}
\usepackage{thmtools,thm-restate}
\newtheorem*{mresult}{Main Result}
\newtheorem{theorem}{Theorem}[section]
\newtheorem{claim}[theorem]{Claim}
\newtheorem{corollary}[theorem]{Corollary}
\newtheorem{lemma}[theorem]{Lemma}
\newtheorem{proposition}[theorem]{Proposition}
\newtheorem{definition}[theorem]{Definition}
\newtheorem{conjecture}[theorem]{Conjecture}
\newtheorem{remark}[theorem]{Remark}

\usepackage[
backend=bibtex,
style=alphabetic,
giveninits=true
]{biblatex}
\DeclareNameAlias{default}{family-given}
\newcommand{\reg}{\mathrm{Reg}}
\newcommand{\ff}{\mathrm{GF}(2)}
\newcommand{\symdiff}{\mathbin{\triangle}}
\newcommand{\bigsymdiff}{\mathop{\vcenter{\hbox{\scalebox{1.55}{$\triangle$}}}}}

\newenvironment{croof}[1][Proof of Claim]{%
  \begin{proof}[#1]%
}{%
  \end{proof}%
}

\DeclareMathOperator{\supp}{supp}

\title{Asymptotically Optimal Bounds for Counting Near-Minimum Circuits in Regular Matroids}
\author{Robert P. Streit\thanks{Department of Electrical and Computer Engineering at The University of Texas at Austin.}\footnotemark[1]\textsuperscript{\;\:,}\thanks{\texttt{rpstreit@utexas.edu}}}
\date{\today}

\begin{document}

\maketitle

\begin{abstract}
    For a regular matroid on $m$ elements with positive weights and weighted girth $g_w$, we prove that the number of circuits of weight at most $\alpha g_w$ is bounded by $(Km)^{2\alpha}$, for every real $\alpha\ge1$ and an absolute constant $K$.
    This improves the previous $m^{O(\alpha^2)}$ bound of Gurjar and Vishnoi and matches classic graph cut and cycle counting bounds of Karger and Subramanian.
    We also show this is asymptotically optimal for all real $\alpha\geq 1$ under the Erdős girth conjecture, while unconditional optimality holds for half-integer $\alpha$ by prior results. 
    As consequences, we resolve a recent conjecture by Khanna, Putterman, and Song implying the existence of a parallel basis search algorithm for regular matroids terminating in logarithmically many rounds of independence queries, and obtain improved bounds on the number of
    near-minimum weight codewords in regular codes and near-shortest vectors in totally unimodular lattices under $\ell_p$-norms.
\end{abstract}

\section{Introduction}

A \emph{circuit} in a matroid is a minimal dependent set, and the \emph{girth} is the size of a smallest circuit.
We ask how the number of circuits grows with their size or weight.
This question was initiated in some sense with Karger's seminal work \cite{karger1993global} showing that, on a graph with $m$ edges, the number of cuts of weight at most $\alpha$ times the minimum is $O(m^{2\alpha})$.
Shortly after, Subramanian proved the analogous statement for graph cycles \cite{subramanian1995polynomial}.
Together, these results show that a graphic or cographic matroid on $m$ elements has $O(m^{2\alpha})$ circuits of size at most $\alpha$ times its girth.
We call these circuits $\alpha$-minimum.
Such results cannot hold for a matroid in general.
For example, there exist binary matroids for which there are exponentially many circuits at the girth \cite{ABV01}.
Counting \(\alpha\)-minimum circuits is fundamental and connects to lattices and error-correcting codes \cite{GV} and to parallel algorithms for combinatorial optimization \cite{khanna_et_al-graphic,gurjar2021isolating}.

By Seymour's theorem \cite{seymour1980decomposition}, every regular matroid can be decomposed into graphic and cographic parts, along with other parts of bounded size.
Using this, Gurjar and Vishnoi \cite{GV} extended the graph bounds to an \(m^{O(\alpha^2)}\) bound for \(\alpha\)-minimum circuits in \(m\)-element regular matroids.
In this work we refine their methods to achieve a bound of $O(m^{2\alpha})$.
So although regular matroids admit a much richer dependency structure than graphs, the asymptotic growth of the number of $\alpha$-minimum circuits behaves no differently than that of graph cuts and cycles.

We analyze the weighted form of this problem.
For a matroid $(E,\mathcal M)$ and weights $w:E\to\mathbb{R}_{>0}$, let $\mathscr{C}(\mathcal M)$ be the circuits and $g_w\triangleq\min_{C\in\mathscr{C}(\mathcal M)}\sum_{y\in C} w(y)$ the weighted girth.
With $\reg$ as the set of all regular matroids, define the extremal counting function
\[
    f_\reg(m,\alpha) \triangleq \max_{\substack{(E,\mathcal M) \in \reg : |E| \leq m,\\ w:E\to\mathbb{R}_{>0}}}\left|\left\{C\in\mathscr{C}(\mathcal{M})\mid w(C) \leq \alpha g_w(\mathcal{M})\right\}\right|,
\]
where we abuse notation to let $w(C) = \sum_{y\in C}w(y)$.
Our main result is as follows. 
\begin{mresult}\label{mr}
    There exists an absolute constant $K\geq 1$ such that $f_\reg(m,\alpha)\leq (Km)^{2\alpha}$
    for every integer $m\geq 1$ and real $\alpha \geq 1$. 
\end{mresult}

Like \cite{GV}, we note that our technique applies not only to regular matroids but to any class of binary matroids decomposing into 1- and 2-sums of regular and bounded-size matroids.
For example, our main result applies to max-flow min-cut matroids \cite{seymour1980decomposition,seymour1977matroids} and binary Rayleigh matroids~\cite{choe2006rayleigh} in particular.
This is elaborated on in Remark~\ref{remark:ext}.

The exponent \(2\alpha\) cannot be improved in general.
Specifically, the complete bipartite graph $K_{n,n}$ on $2n$ vertices has $m=n^2$ edges and girth $4$.
Lemma 3.1 of \cite{durocher2015cycle} shows there are $\Theta(n^{2k})$ cycles of length $2k$ in such graphs when $k$ is fixed.
Hence, by letting $k=2\alpha$ we obtain graphic matroids with $\Theta(m^{2\alpha})$ many $\alpha$-minimum circuits for half-integer $\alpha \geq 1$.
In Section~\ref{sec:lb}, we apply Simonovits' Even-Cycle Supersaturation Theorem \cite{jiang2020supersaturation,erdos1984cube} with the Erdős Girth Conjecture \cite{erdos1964extremal} to conditionally extend this lower bound to real $\alpha$ (see Theorem~\ref{thm:lower-bound}).
So, under this conjecture, our main result obtains the exact polynomial rate of growth for any fixed real $\alpha\geq 1$.

\begin{corollary}\label{cor:tight}
    Assuming the Erdős girth conjecture, for every fixed real $\alpha\geq1$,
    \[
        \lim_{m\to\infty}\frac{\log f_\reg(m,\alpha)}{\log m}=2\alpha.
    \]
    Moreover, this holds for all $\alpha\in\frac 1 2\mathbb{Z}$ with $\alpha\geq 1$ unconditionally.
\end{corollary}

\paragraph{Implications.}

Very recently, Khanna, Putterman, and Song \cite{khanna_et_al-graphic} showed that a bound of $m^{O(\alpha)}$ on the number of $\alpha$-minimum circuits in every deletion minor (relative to its own girth) of a binary matroid suffices to find a basis in $O(\log m)$ adaptive rounds of polynomially many independence queries. 
They conjectured regular matroids satisfy this property (as well as max-flow min-cut matroids). 
Our main result resolves their conjecture with the tighter count $O(m^{2\alpha})$.
Furthermore, it is easy to observe that binary matroids described as 1- and 2-sums of regular and bounded-size matroids are closed under deletions\footnote{Deletion can be performed in the summand containing the deleted element. If this makes a shared element of a 2-sum a coloop, no circuit can cross that sum, which consequently becomes a 1-sum of the two summands with the shared element deleted. The statement follows since regular matroids and bounded-size matroids are minor-closed.}.
So this consequence extends through this class.

\begin{corollary}\label{cor:par}
    Every binary matroid on $m$ elements which decomposes as 1- and 2-sums of regular and bounded-size matroids admits a parallel algorithm outputting a basis, with high probability, in $O(\log m)$ adaptive rounds of polynomially many independence queries.
\end{corollary}


In Section~\ref{sec:app} we also revisit two algebraic applications of $\alpha$-minimum circuit counting identified by Gurjar and Vishnoi \cite{GV}.
The first concerns error-correcting codes.
Specifically, a \emph{regular code} is a binary linear code with a parity-check matrix representing a regular matroid \cite{kashyap2008decomposition}.
Every nonzero codeword of a regular code can be expressed as a disjoint union of circuits of an associated regular matroid, with their sizes summing to the codeword's Hamming weight. 
So, bounds on the number of $\alpha$-minimum circuits constrain the number of codewords whose weight is close to the minimum distance.
For a regular code with distance $d$, \cite{GV} bounds the number of codewords with Hamming weight at most $\alpha d$ by $m^{O(\alpha^3)}$.
We obtain the following improvement.
\begin{theorem}\label{thm:codes}
    For any regular code $\mathcal C\subseteq \ff^m$ of distance $d$ and fixed real $\alpha\geq 1$, the number of codewords with Hamming weight at most $\alpha d$ is no greater than $m^{2\alpha + o(1)}$.
\end{theorem}

The second application concerns lattices defined by the nullspace of totally unimodular matrices.
Specifically, given an $n\times m$ totally unimodular matrix $A$, we call a lattice $\mathcal L(A)\triangleq \ker(A)\cap \mathbb{Z}^m$ a \emph{totally unimodular lattice}.
Every vector in such a lattice can be decomposed into a conformal sum of vectors corresponding to the circuits of the represented matroid $M[A]$.
Because of this, one may control the number of nearly shortest lattice vectors using a bound on the number of $\alpha$-minimum circuits.
Gurjar and Vishnoi \cite{GV} establish an $m^{O(\alpha^6)}$ bound for the Euclidean norm, with their argument extending to $m^{O(\alpha^{3p})}$ for every fixed finite $p\ge1$.
We improve this to $m^{2\alpha^p+o(1)}$.
\begin{theorem}\label{thm:lattices}
    Let $A$ be an $n\times m$ totally unimodular matrix with $\mathcal L(A)\neq\{\mathbf{0}\}$, and define 
    \[
        \lambda_p=\min_{z\in \mathcal L(A)\setminus\{\mathbf 0\}}\|z\|_p.
    \]
    For every fixed $p\geq 1$ and fixed real $\alpha \geq 0$, the number of nonzero vectors of $\mathcal L(A)$ of $\ell_p$-norm at most $\alpha\lambda_p$ is no greater than $m^{2\alpha^p+o(1)}$.
\end{theorem}
These improvements do not follow directly by substituting our main result into \cite{GV}, which incurs additional losses when counting how circuits combine into the objects of interest.
We employ an elementary lemma, Lemma~\ref{lem:ele-cnt}, to avoid these losses.
Taken together, these results reveal that assembling circuits into codewords or lattice vectors introduces no further increase in the polynomial growth rate, after accounting for the natural rescaling by $\alpha^p$ in the $\ell_p$ setting. Thus, the extremal counting behavior of graph cycles and cuts persists through both the combinatorial structure of regular matroids and the algebraic objects they define.

\paragraph{Overview}
Following \cite{GV}, Section~\ref{sec:mr} first selects parts of the tree to direct the counting in Lemma~\ref{lem:C-assignment}, then constructs simpler graph counting subproblems encoding pieces of circuits corresponding to those parts, called \emph{projections}, in Lemma~\ref{lem:mapping}, and finally combines local counting estimates from Lemmas~\ref{lem:other-bound} and~\ref{lem:heavy-D-bound} to obtain the result.
Constructions such as boundaries and signatures are taken from \cite{GV}.
However, important differences arise from how we construct and aggregate our subproblems.
Specifically, \cite{GV} uses their subproblems to uniformly bound the number of configurations of each selected projection.
They then multiply these bounds and incur an $\alpha^2$ in their exponent.
We take a more quantitative approach.
The classical bounds suggest that the number of possibilities in each projection (as controlled by a graph counting subproblem we create) should correspond to its contribution to its circuit's weight, as the polynomial rate of growth of $\alpha$-minimum cuts and cycles is proportional to the weight threshold parameterized by $\alpha$.
So, we select the distinguished projections of a circuit and encode them as cycles and cuts in graphic subproblems with local weights that sum to the weight of the original circuit.
This is given in Lemmas~\ref{lem:C-assignment}.3 and~\ref{lem:mapping}.1.
Combining this idea with an application of a new weight-sensitive counting estimate of cuts and cycles satisfying certain conditions, given in Lemma~\ref{lem:heavy-D-bound}, along with careful analysis gives the main result.

One key difficulty is that the projections we distinguish in Sections~\ref{subsec:assigning} and~\ref{subsec:reduction} may be arbitrarily lighter than the original matroid's weighted girth, and so the classic $\alpha$-minimum cycle and cut bounds will not apply to our subproblems at the right scale.
Gurjar and Vishnoi address this by introducing additional weight to inflate the girth of their graphic subproblems so that \cite{karger1993global,subramanian1995polynomial} can be easily applied.
Yet this is incompatible with our chosen approach, as the additional weight distorts the contribution of each projection through its local weights to its circuit.
To solve this, we present a specialized graph cut and cycle counting bound in Lemma~\ref{lem:heavy-D-bound} showing a lower bound on the weights of cycles and cuts avoiding a prescribed set of edges $R$ is sufficient to control the weighted contributions of those containing $R$.
In Lemma~\ref{lem:mapping}.3, we prove that our subproblems each satisfy this premise so that we may apply Lemma~\ref{lem:heavy-D-bound} in the proof of our main result.
A further challenge is accounting for the different ways selected projections can be chosen, which could introduce additional polynomial factors in the final count.
Lemma~\ref{lem:C-assignment}.2 ensures each selected piece accounts for at least $g_w/2$ of the circuit's weight, which enables an algebraic comparison between the local and global counting expressions compensating for the number of possible selections (see~\eqref{eq:ll} and Claim~\ref{claim:1}).

After finishing the proof of the main result, in Section~\ref{sec:lb} we show the existence of a family of graphic matroids with at least $m^{2\alpha-\epsilon}$ many $\alpha$-minimum circuits, for all $\varepsilon > 0$, under the Erdős girth conjecture in Theorem~\ref{thm:lower-bound}, which gives Corollary~\ref{cor:tight}.
Finally, we conclude our work by proving Theorems~\ref{thm:codes} and~\ref{thm:lattices} in Section~\ref{sec:app} through a shared counting principle given in Lemma~\ref{lem:ele-cnt}.

\paragraph{Related Work}

Classical results of Karger \cite{karger1993global} and Subramanian \cite{subramanian1995polynomial} bound the number of near-minimum cuts and cycles, respectively, with the shortest-cycle case also studied by Teo and Koh \cite{teo1992number}.
Karger's bounds have been an important ingredient in enumeration \cite{KS96} and sampling algorithms \cite{Kar99} for cuts and flows, as well as in network reliability problems \cite{Kar01}.
Seymour's decomposition theorem \cite{seymour1980decomposition}, which led to the first efficient algorithm for testing total unimodularity, has been quite useful for reducing facts about regular matroids to analysis on graphs.
This has been employed for finding minimum-cycle bases \cite{golynski2002polynomial}, approximating the partition function of the ferromagnetic Ising model \cite{goldberg2013polynomial}, the matroid secretary problem \cite{dinitz2014matroid}, obtaining small extended formulations of the independence polytope \cite{aprile2022regular}, and derandomizing the isolation lemma \cite{gurjar2021isolating} to name a few.
Importantly for our work, \cite{gurjar2021isolating} obtained a polynomial bound on the number of circuits of weight strictly less than $3/2$ times the minimum circuit weight, which motivated the more detailed study in \cite{GV}.

Concurrent work \cite{chandrasekaran2026random} independently bounds the $\alpha$-minimum circuits in regular matroids by $m^{O(\alpha)}$, and so also resolves the conjecture of \cite{khanna_et_al-graphic} and obtains Corollary~\ref{cor:par}.
They additionally obtain parallel algorithms for computing the smallest circuit.
Our arguments share some technical ingredients 
but differ fundamentally in their approaches.
Their proof follows the philosophy of Karger's contraction algorithm \cite{karger1993global}, bounding the number of $\alpha$-minimum circuits by analyzing their survival under random deletion and the structure of the resulting matroid.
In contrast, we pursue a quantitative refinement of the direct counting method developed by Gurjar and Vishnoi which carefully accounts for the weights of projections and their combinations to obtain the optimal exponent.

\section{Preliminaries}\label{sec:preliminaries}

For background and terminology surrounding matroids we refer the reader to \cite{oxley2006matroid,welsh2010matroid,white1987combinatorial}.
For a binary matroid $\mathcal M$ let $\mathscr C(\mathcal M)$ denote its circuits and $\mathscr Z(\mathcal M)$ its cycles,
meaning disjoint unions of circuits, including the empty set.
A matroid on a fixed ground set is determined by its cycles.
The cycles $\mathscr Z(\mathcal M)$ are closed under symmetric difference for binary $\mathcal M$ \cite{white1987combinatorial}.
For a graph $G$ we specify the edges and vertices by $\mathcal E(G)$ and $\mathcal V(G)$.
Graphs may contain loops and parallel edges unless stated otherwise. 
In a graphic matroid, circuits correspond to graph cycles.
In a cographic matroid, circuits correspond to minimal nonempty cuts, whereas matroid cycles correspond to general cuts. 
When a decomposition part is both graphic and cographic, we treat it as graphic.

Let $\mathcal M_1$ and $\mathcal M_2$ be binary matroids with nonempty ground sets $E_1$ and $E_2$.
The sum $\mathcal N=\mathcal M_1\symdiff\mathcal M_2$ is another binary matroid over the ground set $E_1\symdiff E_2$ whose cycle family is
\begin{equation}\label{eq:ksum-def}
    \mathscr Z(\mathcal N) = \{C_1\symdiff C_2\subseteq E_1\symdiff E_2\mid (C_1,C_2)\in \mathscr Z(\mathcal M_1)\times\mathscr Z(\mathcal M_2)\}.
\end{equation}
The intersection $E_1\cap E_2$ is called the \emph{separator}.
The separator elements cancel in~\eqref{eq:ksum-def} and do not belong to the ground set of \(\mathcal N\); they serve only to specify the sum.
Special cases of this operation are called $k$-sums.
Let the separator $S\triangleq E_1\cap E_2$.
When $S=\varnothing$, then we write $\mathcal N=\mathcal M_1\oplus_1\mathcal M_2$ and say $\mathcal N$ is the 1-sum of $\mathcal M_1$ and $\mathcal M_2$.
If $|S|=1$ and its element is not a loop or coloop of $\mathcal M_1$ or $\mathcal M_2$ (and $|E_1|,|E_2|\ge3$) then $\mathcal N=\mathcal M_1\oplus_2\mathcal M_2$ is a 2-sum.
Finally, if $|S|=3$ and $S$ is a circuit of both summands while containing no cocircuit of either (and $|E_1|,|E_2|\ge7$), then $\mathcal N=\mathcal M_1\oplus_3\mathcal M_2$ is a 3-sum.
A matroid is \emph{3-connected} if it is not a 1- or 2-sum.

Seymour’s theorem \cite{seymour1980decomposition} expresses every regular matroid as repeated 1-, 2-, and 3-sums of graphic, cographic, and $\mathcal R_{10}$ matroids.
Following Gurjar and Vishnoi \cite{GV}, we employ a decomposition tree that permits these sums to be taken in arbitrary orders.
\begin{definition}[Unordered Decomposition Tree]\label{def:udt}
    An unordered decomposition tree (UDT) of a binary matroid \((E,\mathcal M)\) is a tree \(T\) with a binary matroid \((E_v, M_v)\), $E_v\neq\varnothing$, associated with each \(v\in\mathcal V(T)\), satisfying the following conditions.
    \begin{enumerate}
        \item $|E_u\cap E_v| \in \{0,1,3\}$ for all $(u,v)\in \mathcal E(T)$.
        \item $E_u\cap E_v = \varnothing$ for all $(u,v) \notin \mathcal E(T)$ with $u\neq v$.
        \item For every nonempty \(X\subseteq\mathcal V(T)\) inducing a subtree, there is a unique matroid \((E_X,\mathcal M_X)\), where \(E_X=\mathop{\triangle}_{u\in X}E_u\) and \(\mathcal M_{\{u\}}=\mathcal M_u\), such that $\mathcal M_X=\mathcal M_{X_1}\oplus_k\mathcal M_{X_2}$
for any vertex sets $X_1$ and $X_2$ of components obtained by deleting an edge of \(T[X]\). Moreover, $(E,\mathcal M)=\left(E_{\mathcal V(T)}, M_{\mathcal V(T)}\right).$
    \end{enumerate}
\end{definition}

For an edge $(u,v)\in\mathcal E(T)$, we let $S_{u,v}$ be the separator in the corresponding $k$-sum of $M_u=M_{\{u\}}$ and $M_v=M_{\{v\}}$.
For the vertices $X$ corresponding to a subtree, we let \(N(X)\) denote the vertices outside \(X\) adjacent to a vertex of \(X\), and $S_{u,X}$ give the separator for the sum of $M_X$ and $M_u$ when $u\in N(X)$.
Gurjar and Vishnoi obtain a UDT using \cite{chun2011unavoidable}.
We use a refinement of 3-connected parts from \cite{aprile2022regular} to obtain cographic structure needed in Lemma~\ref{lem:mapping}. 
The proof is deferred to Appendix~\ref{app:decomp}.

\begin{restatable}{proposition}{restatefoundation}\label{prop:foundation}
Every regular matroid $\mathcal M$ on $m\geq 1$ elements
admits a UDT $T$ whose vertices correspond to graphic,
cographic, or $\mathcal R_{10}$ matroids, and whose edges
correspond to $1$-, $2$-, or $3$-sums.
Moreover, $T$ may be chosen so that:
\begin{enumerate}
    \item For every $u\in\mathcal V(T)$ where $(E_u,M_u)$ is cographic, there exists a
    graph $G_u$ with $\mathcal E(G_u) = E_u$ representing
    $\mathcal M_u$ cographically, such that every
    three-element separator $S_{u,v}$ incident to $u$
    is the set of edges incident to a degree-three
    vertex of $G_u$.
    \item The decomposition satisfies $\sum_{u\in\mathcal V(T)}|E_u|\leq 7m-6.$
    \item Every vertex $u$ with $ M_u\cong\mathcal R_{10}$ is incident only to $1$- and $2$-sum edges.
\end{enumerate}
\end{restatable}


The following collection of observations and results from \cite{GV} expresses each circuit as a symmetric difference of circuits from the parts of the decomposition or empty sets, called its projections.
%

\begin{proposition}[\cite{GV}]\label{prop:projections}
    Let $T$ be a UDT of a binary matroid $(E,\mathcal M)$ with weights $w:E\to\mathbb{R}$.
    Every circuit $C$ of $\mathcal M$ can be associated with projections $(C_u)_{u \in \mathcal V(T)}$ where each $C_u$ is either empty or a circuit of $M_u$, such that we may write $C = \mathop{\triangle}_{u \in \mathcal V(T)} C_u$.
    For $X\subseteq\mathcal V(T)$, define $C_X\triangleq\mathop{\triangle}_{u\in X} C_u$.
    \begin{enumerate}
        \item The support $\{u \in \mathcal V(T)\mid C_u \neq \varnothing\}$ forms a connected subtree $T'$.
        \item For all $(u,v)\in\mathcal E(T)$ we have that $C_u\cap S_{u,v} = C_v\cap S_{u,v}$. The common intersection is empty or a singleton, and is nonempty whenever both \(C_u\) and \(C_v\) are nonempty.
        \item The sets \(C_u\cap E\) are pairwise disjoint, have union \(C\), and $w(C_X\cap E) = \sum_{u \in X} w(C_u \cap E)$ for all $X\subseteq\mathcal V(T).$ 
    \end{enumerate}
\end{proposition}
\section{Counting Near-Minimum Weighted Circuits}\label{sec:mr}
%
This section is devoted to the proof of the main result.
In what follows, fix a regular matroid $(E,\mathcal M)$, weights $w:E\to\mathbb R_{> 0}$, and a UDT $T$ with the structure given in Proposition~\ref{prop:foundation}.

\subsection{Assigning Circuits to Pieces of the UDT}\label{subsec:assigning}

For each circuit \(C\), we select a set of vertices \(\mathcal S(C)\) according to the weights of its projections. These vertices index the local graph-counting problems constructed in Lemma~\ref{lem:mapping}.

Let $X\subseteq\mathcal V(T)$ be a set of vertices defining a connected subtree of $T$.
It will be important to record how a circuit's projection over $X$ interacts with the rest of the UDT.
For any cycle \(Z\in\mathcal Z(\mathcal M_X)\), define its \emph{boundary} by
\[
    \partial_X Z \triangleq \bigcup_{u\in N(X)}Z \cap S_{(u,X)}.
\]
The boundary thus records how \(Z\) meets the rest of the UDT.
%
The following generalization of Claim 4.5 from \cite{GV} shows that each boundary determines at most one cycle whose weight on the original elements is below \(g_w/2\).
Essentially, light portions of a circuit are determined by how they connect to the rest of the decomposition, allowing us to recover them without counting their possibilities separately.
\begin{lemma}\label{lem:bd=>equal}
    Let $T$ be a UDT of $\mathcal M$ and $X\subseteq\mathcal V(T)$ induce a connected subtree.
    Suppose $Z_1,Z_2\in\mathscr Z(M_X)$ are two cycles of $M_X$ such that $\partial_X Z_1 = \partial_X Z_2$ and $\max_{i\in\{1,2\}}w(Z_i\cap E) < g_w/2$.
    Then $Z_1 = Z_2$.
\end{lemma}
\begin{proof}
    By way of contradiction suppose $Z_1\neq Z_2$.
    Because $M_X$ is a binary matroid, $Z_1\symdiff Z_2$ is a nonempty cycle of $M_X$.
    Observe that
    \begin{equation*}
        \begin{split}
            \partial_X(Z_1\symdiff Z_2) &= \bigcup_{u\in N(X)}(Z_1\symdiff Z_2)\cap S_{u,X}\\
                                                  &= (Z_1\symdiff Z_2)\cap \left(\bigcup_{u\in N(X)}S_{u,X}\right)\\
                                                  &= \left(Z_1\cap \bigcup_{u\in N(X)}S_{u,X}\right)\symdiff\left(Z_2\cap \bigcup_{u\in N(X)}S_{u,X}\right)\\
                                                  &=\partial_X Z_1 \symdiff \partial_X Z_2.
        \end{split}
    \end{equation*}
    So equality of the boundaries implies $\partial_X(Z_1\symdiff Z_2)=\varnothing$.
    Therefore $Z_1\symdiff Z_2$ contains no separator elements appearing in ground sets of matroids in the decomposition outside of $(E_X,M_X)$, and so $Z_1\symdiff Z_2$ is a cycle of $\mathcal M$.
    Consequently,
    \[
        g_w \leq w(Z_1\symdiff Z_2) \leq w(Z_1\cap E) + w(Z_2 \cap E) < g_w,
    \]
    a contradiction.
\end{proof}

Root \(T\) by a fixed rule independent of \(C\), and process its vertices in postorder. At \(u\), let \(P_u\) be the union of \(\{u\}\) and the sets returned by its children.
The selection rule is as follows:
if \(w(C_{P_u}\cap E)\ge g_w/2\), add \(u\) to $\mathcal S(C)$ and return the empty set to its parent.
Otherwise, when $w(C_{P_u}\cap E) < g_w/2$, leave $\mathcal S(C)$ unchanged and return the set $P_u$.
As a remark, since $w(C_{\mathcal V(T)}\cap E) = w(C)\geq g_w$ it follows that $|\mathcal S(C)| \geq 1$ upon termination.
A similar construction is given in \cite{GV} using a centroid decomposition to obtain a set $U$ such that every component $X$ of $T-U$ satisfies $w(C_X\cap E)< g_w/2$, so that an analogue of Lemma~\ref{lem:bd=>equal} may be applied to validate that one need only concern themselves with how the circuit projects as $C_u$ onto each $u\in\mathcal S(C)$ for the purposes of counting.
We prove this property is maintained by our method in Lemma~\ref{lem:C-assignment}.1, along with some additional properties enabling the ideas in our overview.

Call the components of \(T-\mathcal S(C)\) \emph{unselected}. Assign each component \(X\) to an adjacent selected vertex \(a_{\mathcal S(C)}(X)\), whose local counting problem will account for \(X\)’s weight contribution.
Select $a_{\mathcal S(C)}(X)$ in the following manner: If $X$ does not contain the root, then $a_{\mathcal S(C)}(X)$ is the parent of the vertex in $X$ closest to the root.
Otherwise, when $X$ contains the root, $a_{\mathcal S(C)}(X)$ is a neighbor $u\in N(X)$ chosen by a fixed tie-breaking rule.
For \(u\in\mathcal S(C)\), define
\[
    B_{u}(C)\triangleq \{u\}\cup\bigcup_{X\in a_{\mathcal S(C)}^{-1}(u)}X.
\]
Lemma~\ref{lem:C-assignment}.3 shows summing $w(C_{B_u(C)}\cap E)$ over $u\in\mathcal S(C)$ gives exactly $w(C)$, while Lemma~\ref{lem:C-assignment}.2 maintains that each weight $w(C_{B_u(C)}\cap E)$ itself is at least $g_w/2$.

\begin{lemma}\label{lem:C-assignment}
    Let $T$ be a UDT of $\mathcal M$.
    For all circuits $C \in \mathscr{C}$, the construction has the following properties:
    \begin{enumerate}
        \item $w(C_X\cap E) < g_w/2$ for each component $X$ of $T-\mathcal S(C)$.
        \item $w(C_{B_u(C)}\cap E) \geq g_w/2$ and $C_u \neq \varnothing$ for all selected $u\in\mathcal S(C)$.
        \item $\sum_{u \in \mathcal S(C)} w(C_{B_u(C)}\cap E) = w(C)$.
    \end{enumerate}
\end{lemma}
\begin{proof}
    First observe that the sets $P_u$ formed at selected vertices $u\in\mathcal S(C)$ are pairwise disjoint, and all vertices not contained in any such set either form the empty set or a component containing the root.
    Call this latter set $Y$; we have $w(C_Y\cap E)<g_w/2$ if $Y\neq \varnothing$ as otherwise a vertex of $Y$ would have been selected.
    For 1. select a component $X$ in $T-\mathcal S(C)$.
    If $X$ contains the root then $X=Y$ and our prior observation shows $w(C_X\cap E) < g_w/2$.
    Otherwise, if $X$ does not contain the root, then examine the vertex $v\in X$ closest to the root.
    Observe that when \(v\) was visited in the postorder traversal, it returned the set $X$ to its parent.
    This is because $X$ contains exactly those unselected descendants connected to $v$, since all selected descendants had returned empty sets.
    But then for $v\notin\mathcal S(C)$ to hold it must be that $w(C_X\cap E) < g_w/2$.

    Fix \(u\in\mathcal S(C)\). Each component returned by a child of \(u\) is assigned to \(u\), so \(P_u\subseteq B_u(C)\). Hence,
    \begin{equation}\label{eq:bupugw2}
        w(C_{B_u(C)}\cap E)\geq w(C_{P_u}\cap E) \geq g_w/2,
    \end{equation}
    where the last inequality follows from the fact that $u$ was selected.
    Since every unselected component is assigned to exactly one selected vertex, the sets $\{B_u(C)\mid u\in\mathcal S(C)\}$ partition $\mathcal V(T)$.
    Therefore, Proposition~\ref{prop:projections}.3 gives
    \[
        \sum_{u\in\mathcal S(C)}w(C_{B_u(C)}\cap E) = \sum_{u\in\mathcal V(T)} w(C_u\cap E) = w(C).
    \]
    This proves 3. and the first part of 2.

    What remains for 2. is showing $C_u\neq \varnothing$ for all $u\in\mathcal S(C)$.
    Suppose $u$ is selected and $C_u=\varnothing$.
    Each set returned by a child has weight strictly less than $g_w/2$.
    Since their union with $\{u\}$ has weight at least $g_w/2$, at least two distinct children returned sets of positive weight.
    Choose a vertex $x_i$ with $w(C_{x_i}\cap E)>0$ in each of two such sets.
    The path from $x_1$ to $x_2$ passes through $u$.
    By Proposition~\ref{prop:projections}.1 the set of vertices with nonempty projections forms a connected subtree, so we have the contradiction $C_u\neq\varnothing$.
%
%
\end{proof}

\subsection{Reduction to Weighted Cycle and Cut Counting}\label{subsec:reduction}

Lemmas~\ref{lem:bd=>equal} and~\ref{lem:C-assignment}.1 reduce the task to counting projections at selected vertices. The \(M_u\cong\mathcal R_{10}\) case has bounded size.
For graphic and cographic pieces, we encode projections as cycles and cuts in weighted graphs shared by suitable families of circuits.
Grouping circuits only by their selected vertices is insufficient, since the boundaries of projections onto unselected components adjacent to multiple selected vertices depend on how a circuit connects them.
We instead group by a \emph{signature} recording these connections, thereby fixing the corresponding boundaries.
Then, Lemmas~\ref{lem:bd=>equal} and~\ref{lem:C-assignment}.1 apply for unselected components adjacent to multiple selected vertices, while the remaining unselected components adjacent to a single selected vertex are handled separately.
Moreover, via Lemma~\ref{lem:C-assignment}.3 we can assign weights in the graphs so that their local objects (cuts and cycles) can be used to recover each circuit's weight (see Lemma~\ref{lem:mapping}.1).
Finally, to apply the local estimates dictated by Lemma~\ref{lem:heavy-D-bound}, we must ensure that every local cycle or nonempty cut avoiding connections to other selected vertices extends to a nonempty cycle of the original matroid with the same weight, and hence has weight at least $g_w$. 
Our construction guarantees this property (see Lemma~\ref{lem:mapping}.3).

Fix a circuit $C\in\mathcal C(\mathcal M)$, and for a selected vertex $u\in \mathcal S(C)$ let $ N_u^*(\mathcal S(C))$ consist of neighbors of $u$ leading towards another selected vertex,
\[
     N_u^*(\mathcal S(C))\triangleq\{v \in N(u)\mid (\exists w \in \mathcal S(C)\setminus\{u\})\;v\text{ lies on some path from $u$ to $w$}\}.
\]
Recall that the vertices $v\in\mathcal V(T)$ with a nonempty projection $C_v\neq \varnothing$ form a subtree by Proposition~\ref{prop:projections}.1.
Since Lemma~\ref{lem:C-assignment}.2 gives $C_u\neq \varnothing$, this implies that any vertex on a path between $u$ and another selected vertex must also admit a nonempty projection, i.e.,
\begin{equation}\label{eq:neighbors-cv-nonempty}
    u \in \mathcal S(C)\text{ and }v \in  N_u^*(\mathcal S(C)) \implies C_v \neq \varnothing.
\end{equation}
Then, use $ N_u^*(\mathcal S(C))$ to collect separator elements appearing in $C_u$ from edges in $T$ leading towards other selected vertices by defining $R_u(C) \triangleq \bigcup_{v \in  N_u^*(\mathcal S(C))}(C_u \cap S_{u,v}).$
Importantly, since $u$ is a neighbor to any $v \in  N_u^*(\mathcal S(C))$, applying \eqref{eq:neighbors-cv-nonempty} to Proposition~\ref{prop:projections}.2 implies
\begin{equation}\label{eq:proj-sep-single}
    u \in \mathcal S(C)\text{ and }v \in  N_u^*(\mathcal S(C)) \implies |C_u\cap S_{u,v}| = 1.
\end{equation}
In particular, $R_u(C)$ is empty if and only if $ N_u^*(\mathcal S(C))$ is empty.
With that, the \emph{signature} of a circuit $C$ is the collection of tuples $\sigma(C)\triangleq\{(u,R_u(C))\mid u\in\mathcal S(C)\}$.
For a signature $A\in\sigma(\mathscr C(\mathcal M))$, call its first coordinates $\mathcal S(A)$ and let $R_u(A)$ be the value paired with $u \in\mathcal S(A)$.
Our construction follows.

\begin{lemma}\label{lem:mapping}
    For every signature \(A\in\sigma(\mathscr C(\mathcal M))\), there exist an injection
    \[\sigma^{-1}(A)\to \bigtimes_{u\in \mathcal S(A)}2^{E_u}:C\mapsto (D_u^{(A)}(C))_{u\in \mathcal S(A)},\]
    weight functions $\omega_u^{(A)}:E_u\to\mathbb{R}_{\geq 0}$ for $u\in\mathcal S(A)$, and graphs $G_u^{(A)}$ for $u\in \mathcal S(A)$ where $ M_u$ is graphic or cographic satisfying $R_u(A)\subseteq\mathcal E(G_u^{(A)})\subseteq E_u$, such that the following are all true:
    \begin{enumerate}
        \item For all $u \in \mathcal S(A)$ and $C\in\sigma^{-1}(A)$, $R_u(C)\subseteq D_u^{(A)}(C)$ and $\omega_u^{(A)}\!(D_u^{(A)}(C)) = w(C_{B_u(C)}\cap E)$.
        \item For each graphic (resp. cographic) node $u\in \mathcal S(A)$ and $C \in \sigma^{-1}(A)$, $D_u^{(A)}(C)$ is a cycle (resp. nonempty cut) of $G_u^{(A)}$.
        \item For each graphic (resp. cographic) node $u\in \mathcal S(A)$, every cycle (resp. nonempty cut) $D$ of $G_u^{(A)}$ satisfying $D\cap R_u(A)=\varnothing$ has $\omega_u^{(A)}(D)\geq g_w$.
    \end{enumerate}
\end{lemma}
\begin{proof}
    Fix a signature \(A\), and let \(\mathcal X\) be the components of \(T-\mathcal S(A)\). Every \(X\in\mathcal X\) satisfies \(N(X)\subseteq\mathcal S(A)\). We distinguish components with one selected neighbor from those with at least two.
    In the latter case, the signature fixes the boundary.
    \begin{claim}\label{claim:inc2-boundary}
        If $X$ is incident to at least two vertices from $\mathcal S(A)$, then $\partial_X C_X$ is determined by $A$ for any circuit $C$ with signature $A$. 
        To be precise,
        \begin{equation}\label{eq:boundary-fixed-by-A}
            | N(X)|\geq 2\implies (\forall C\in\sigma^{-1}(A))\;\partial_XC_X = \bigcup_{u\in N(X)} R_u(A)\cap S_{u,X}.
        \end{equation}
    \end{claim}
    \begin{croof}
        Select $C\in\sigma^{-1}(A)$.
        If $X$ is adjacent to at least two selected vertices in $\mathcal S(A)$, then any $v\in X$ adjacent to $u\in\mathcal S(A)$ lies on a path toward another selected vertex.
        Hence $v\in N_u^*(\mathcal S(C))$, so $R_u(C)$ includes the set $C_u\cap S_{u,X}$.
        Then because separators are disjoint, we observe $R_u(C) \cap S_{u,X} = C_u\cap S_{u,X}$.
        The definition of a UDT also restricts the separator $S_{u,X}$ to only lie in the ground sets $E_u$ and $E_v$ in the decomposition, so only $C_v$ can contribute elements of $S_{u,X}$ to $C_X$.
        It follows that $C_X\cap S_{u,X}=C_v\cap S_{u,X}$.
        Then by Proposition~\ref{prop:projections}.2, we've identified $C_X\cap S_{u,X}=R_u(C) \cap S_{u,X}$.
        As $R_u(A) = R_u(C)$ we have the claim.
    \end{croof}

    Components $X$ adjacent to a single selected vertex require a different approach.
    We instead show every singleton from a separator $S_{u,X}$ for $u\in\mathcal S(A)$ is the boundary of some cycle of $M_X$.
%
%
    \begin{claim}\label{claim:y-is-boundary}
        If $X$ is incident to one vertex from $\mathcal S(A)$, then each element in the corresponding separator defines a singleton boundary determined by some cycle of $M_X$.
        This is to say,
        \[
            \{u\}= N(X)\implies (\forall y \in S_{u,X})(\exists Z\in\mathscr Z(M_X))\;\partial_XZ=\{y\}.
        \]
    \end{claim}
    \begin{croof}
        Suppose $S_{u,X}$ encodes a 2-sum or a 3-sum (so that it is nonempty) and select $y\in S_{u,X}$.
        By the definition of $k$-sums the separator cannot contain a cocircuit of $M_X$, and so $E_X\setminus S_{u,X}$ spans $M_X$.
        Hence there is a basis $I\subseteq E_X\setminus S_{u,X}$, and so we let $Z$ be the fundamental circuit of $y$ with respect to $I$.
        This makes $Z\cap S_{u,X}=\{y\}$, and as $u$ defines the only edge incident to $X$ we find $\partial_X Z = \{y\}$.
    \end{croof}
    To assign weights that faithfully represent the projections on unselected components, we associate each boundary with the minimum weight of a cycle realizing it. 
    Specifically, for a realizable boundary \(\Gamma\), define
    \begin{equation}\label{eq:muX}
        \mu_X(\Gamma)\triangleq\min_{Z\in\mathscr Z(M_X):\partial_X Z = \Gamma}w(Z\cap E).
    \end{equation}
    These weights are defined for every realizable boundary, including those absent from a given signature class, as needed for property 3. For actual circuit projections, the next claim shows that the minimum in~\eqref{eq:muX} is attained by the projection itself.
    \begin{claim}\label{claim:cminimizer}
        Let $Z\in\mathscr Z(M_X)$ and $C\in\sigma^{-1}(A)$ be such that $w(Z\cap E)=\mu_X(\partial_X C_X)$ and $\partial_XZ=\partial_XC_X$.
        Then, $C_X = Z$.
    \end{claim}
    \begin{croof}
        Since $\mathcal S(A) = \mathcal S(C)$, we first observe that $w(C_X\cap E) < g_w/2$ by way of Lemma~\ref{lem:C-assignment}.1.
        Because $C_X$ is a cycle in $M_X$ it is feasible for the minimum in~\eqref{eq:muX}, so $\mu_X(\partial_X C_X) \leq w(C_X\cap E)$.
        Hence $w(Z\cap E)<g_w/2$.
        Lemma~\ref{lem:bd=>equal} and $\partial_XZ=\partial_X C_X$ give the conclusion.
    \end{croof}
    Now for any $u\in\mathcal S(A)$ we define a new weight function $\omega^{(0)}_u$ using the boundary information as follows:
    \begin{equation}\label{eq:hatpi-def}
        \omega^{(0)}_u(y)\triangleq
        \begin{cases}
            w(y),
            \quad& \text{if }y\in E,\\[4pt]
            \displaystyle\mu_X\left(
                \bigcup_{v\in N(X)}
                (R_v(A)\cap S_{v,X})
            \right),
            \quad& \begin{aligned}[c]
                &\text{if }(\exists X\in\mathcal X)\;
                |N(X)|\geq 2,\;
                a_{\mathcal S(A)}(X)=u,\\
                &\quad\text{and }
                \{y\}=R_u(A)\cap S_{u,X},
              \end{aligned}\\[6pt]
            \mu_X(\{y\}),
            \quad& \text{if }(\exists X\in\mathcal X)\;
                N(X)=\{u\}\text{ and }y\in S_{u,X},\\
            0,
            \quad& \text{otherwise}.
        \end{cases}
    \end{equation}
    Claims~\ref{claim:inc2-boundary} and~\ref{claim:y-is-boundary} show each boundary used in~\eqref{eq:hatpi-def} is realizable.
    The next claim shows \(\omega_u^{(0)}\) preserves the contribution of \(B_u(C)\).
    Our strategy is to verify the weight assigned to each separator element in $C_u$ reproduces the weight $w(C_X\cap E)$ of the unselected component $X$ assigned to $u$.
    \begin{claim}\label{claim:hatpi}
        $\omega^{(0)}_u(C_u)=w(C_{B_u(C)}\cap E)$ for all $C\in\sigma^{-1}(A)$.
    \end{claim}
    \begin{croof}
        Fix a circuit $C\in\sigma^{-1}(A)$.
        Firstly, note that if $y \notin E$ and $y\notin S_{u,X}$ for any component $X\in\mathcal X$ with $a_{\mathcal S(A)}(X) = u$, then $\omega^{(0)}_u(y) = 0$.
        So, noting $\mathcal S(C) = \mathcal S(A)$ by construction,
        \begin{equation}\label{eq:picu-decomp}
            \omega^{(0)}_u(C_u) = \omega^{(0)}_u(C_u\cap E) + \sum_{X\in\mathcal X:a_{\mathcal S(C)}(X) = u}\omega^{(0)}_u(C_u\cap S_{u,X}).
        \end{equation}
        Examine any $X\in\mathcal X$ with $a_{\mathcal S(A)}(X) = u$.
        We claim that
        \begin{equation}\label{eq:om2mu}
            \omega_u^{(0)}(C_u\cap S_{u,X})=\mu_X(\partial_X C_X).
        \end{equation}
        There are two cases.
        First suppose $N(X)=\{u\}$.
        Then $S_{u,X}$ is the only separator on the boundary of $X$, and the projections agree across it. 
        Hence $\partial_X C_X=C_u\cap S_{u,X}$.
        If this intersection is $\{y\}$, the third case of~\eqref{eq:hatpi-def} gives $\omega_u^{(0)}(y)=\mu_X(\{y\})$.
        If it is empty then both sides of~\eqref{eq:om2mu} are zero since the empty cycle attains $\mu_X(\varnothing)=0$.
        Now examine the second case $|N(X)|\ge2$.
        The neighbor of $u$ in $X$ lies on the path from $u$ to another selected vertex, so~\eqref{eq:proj-sep-single} gives
        \[
            C_u\cap S_{u,X}=R_u(A)\cap S_{u,X}=\{y\}
        \]
        for some $y$. Since $X$ is assigned to $u$, the second case of~\eqref{eq:hatpi-def} assigns to this one element the cost of the entire boundary.
        And so,
        \[
        \omega_u^{(0)}(y)
        =\mu_X\!\left(\bigcup_{v\in N(X)}(R_v(A)\cap S_{v,X})\right)
        =\mu_X(\partial_X C_X),
        \]
        where the last equality is given by Claim~\ref{claim:inc2-boundary}. This proves~\eqref{eq:om2mu} for every component assigned to $u$.
        Combining this with Claim~\ref{claim:cminimizer} gives
        \begin{equation*}
            \begin{split}
                \omega^{(0)}_u(C_u\cap E) + \sum_{X\in\mathcal X:a_{\mathcal S(A)}(X) = u}\omega^{(0)}_u(C_u\cap S_{u,X}) &= w(C_u\cap E) + \sum_{X\in\mathcal X:a_{\mathcal S(C)}(X) = u}\mu_X(\partial_X C_X)\\
                                                                                                                      &= w(C_u\cap E) + \sum_{X\in\mathcal X:a_{\mathcal S(C)}(X) = u}w(C_X\cap E)\\
                                                                                                                      &= \sum_{v \in B_u(C)}w(C_v\cap E)\\
                                                                                                                      &= w(C_{B_u(C)}\cap E).
            \end{split}
        \end{equation*}
        As our choice of $C\in\sigma^{-1}(A)$ was arbitrary, applying this to~\eqref{eq:picu-decomp} gives the claim.
    \end{croof}

    Fix a vertex $u\in\mathcal S(A)$ such that $(E_u,M_u)$ is either graphic or cographic and a corresponding representation $G^{(0)}_u$ with $\mathcal E(G^{(0)}_u) = E_u$.
    If $M_u$ is cographic but not graphic then by Proposition~\ref{prop:foundation}.1 we may assume every three-element separator $S_{u,v}$ gives the edges of a star in $G^{(0)}_u$.
    Likewise, if $M_u$ is graphic then such a separator is a triangle in $G^{(0)}_u$ since a graphic representation must maintain circuits as cycles.
    Initially, we weight edges according to $\omega^{(0)}_u$.

    The previous claim accounts for the weight of each circuit projection.
    To achieve 3., we modify $G_u^{(0)}$ so that every local cycle or nonempty cut avoiding $R_u(A)$ lifts to a nonempty cycle of $\mathcal M$ of the same weight.
    Its weight is then at least $g_w$.
    Claim~\ref{claim:y-is-boundary} will give such a construction through unselected components whenever they meet the corresponding separators in at most one element.
    However, arbitrary local cycles and cuts need not satisfy this restriction at three-element separators.
    We therefore perform $\Delta$-$Y$ and $Y$-$\Delta$ exchanges (see e.g. Ch. 4 of \cite{truemper1992matroid}) so that, upon reversing these transformations, every local object meets the relevant separators in at most one element (we remark that the independent concurrent work \cite{chandrasekaran2026random} performs these operations as well for similar reasons).
    Finally, deletion and contraction will eliminate connections inconsistent with the signature, allowing us to lift local objects avoiding $R_u(A)$ without involving other selected vertices.

    Define $\mathcal X_3(u)\triangleq \{X\in\mathcal X\mid  N(X)=\{u\}\text{ and }|S_{u,X}|=3\}$, examine some $X\in\mathcal X_3(u)$, and let $S_{u,X}=\{y_1,y_2,y_3\}$.
    If $M_u$ is graphic, meaning $S_{u,X}$ is a triangle in $G_u^{(0)}$, perform a $\Delta$-$Y$ exchange by deleting the edges labelled by $S_{u,X}$, adding a new vertex $v'_X$, and adding an edge $(v'_X,v_i)$ for each of the three $v_i$ vertices from the original triangle.
    Observe that a cycle meeting $S_{u,X}$ at $\{y_i\}$ in $G^{(0)}_u$ is transformed by replacing $(v_j,v_k)$ with the path $v_j,v'_X,v_k$.
    Finally, reintroduce the label `$y_i$' by reassigning it to the edge $(v'_X,v_i)$ when $v_i$ was opposite $y_i$ in the original triangle.
    If $M_u$ is instead cographic, meaning $S_{u,X}$ is a star in $G_u^{(0)}$, a $Y$-$\Delta$ exchange deletes the center vertex and its incident edges and adds edges $(v_j,v_k)$, where $v_i$ is the vertex initially incident to $y_i$ in the star, to form a triangle.
    Correspondingly any cut meeting the original star in $\{y_i\}$ transforms into a cut meeting the new triangle in $\{(v_i,v_j),(v_i,v_k)\}$.
    Label the new edge $(v_j,v_k)$ by $y_i$.

    Perform these exchanges for every $X\in\mathcal X_3(u)$.
    Distinct separators have disjoint edge sets by the definition of a UDT, so the operations are compatible.
    Call the resulting weighted graph $G^{(1)}_u$, taking note that our relabelling preserves $\mathcal E(G^{(1)}_u) = E_u$.
    Collect the elements from $E_u$ in separators leading towards other selected vertices not appearing in any $C_u$ for $C\in\sigma^{-1}(A)$ in the set
    \[
        F_u\triangleq \left(\bigcup_{v\in N^*_u(\mathcal S(A))} S_{u,v}\right)\setminus R_u(A).
    \]
    The final graph is given by
    \[
        G^{(A)}_u\triangleq \begin{cases}
            G^{(1)}_u - F_u,\quad&\text{if $M_u$ is graphic},\\
            G^{(1)}_u/F_u,\quad&\text{if $M_u$ is cographic}.
        \end{cases}
    \]
    Deletion and contraction retain, respectively, cycles and cuts avoiding \(F_u\).

    The construction depends only on \(A\) and \(u\) and removes no element of \(R_u(A)\).
    %
    Hence $R_u(A) \subseteq \mathcal E(G^{(A)}_u)\subseteq \mathcal E(G^{(1)}_u) = E_u$.
    In both the graphic and cographic cases, a local object meeting a three-element separator $S_{u,X}$ in a single element $y$ is transformed into an object containing the complementary pair $S_{u,X}\setminus\{y\}$, while a projection avoiding the separator remains unchanged.
    We describe this effect by the mapping \(\Phi_u:2^{E_u}\to2^{E_u}\) defined by
    \begin{align}\label{eq:Phi-def}
        \Phi_u(D)\triangleq D\symdiff\left(\bigsymdiff_{X\in\mathcal X_3(u):D\cap S_{u,X}\neq\varnothing} S_{u,X}\right),
    \end{align}
    which will allow us to both perform and reverse these exchanges symbolically.
    We reassign the weights as follows:
    \begin{equation}
        \omega^{(A)}_u(y)\triangleq \begin{cases}
            \frac{\omega^{(0)}_u(S_{u,X})}{2}-\omega^{(0)}_u(y),\quad&\text{if }y\in S_{u,X}\text{ and }X\in\mathcal X_3(u),\\
            0,\quad&\text{if }y\in F_u,\\
            \omega^{(0)}_u(y),\quad&\text{otherwise}.
        \end{cases}
    \end{equation}
    To prove \(\omega_u^{(A)}\) is nonnegative, it suffices to consider elements of \(S_{u,X}\) for \(X\in\mathcal X_3(u)\).
    Let $S_{u,X} = \{y_1, y_2, y_3\}$.
    By Claim~\ref{claim:y-is-boundary}, $\{y_i\}$ is a boundary generated by some circuit over $X$ so it is valid to write and inspect $\mu_X(\{y_i\})$.
    Let $Z_i\in\mathscr Z(M_X)$ be a cycle such that $\partial_X Z_i = \{y_i\}$ and $w(Z_i\cap E) = \mu_X(\{y_i\})$.
    Importantly, this selection makes $Z_i\cap S_{u,X}=\{y_i\}$.
    For the sum $M_u\oplus_3 M_X$ the separator satisfies $S_{u,X}\in\mathscr C(M_u)\cap\mathscr C(M_X)$.
    So,
    we see $Z_j\symdiff Z_k\symdiff S_{u,X}$ is a cycle of $M_X$ with boundary $\{y_i\}$.
    Hence, as $S_{u,X}\cap E=\varnothing$,
    \[
        \mu_X(\{y_i\})\leq w((Z_j\symdiff Z_k\symdiff S_{u,X})\cap E) \leq w(Z_j\cap E) + w(Z_k\cap E) = \mu_X(\{y_j\}) + \mu_X(\{y_k\}),
    \]
    where we let $S_{u,X}\setminus\{y_i\} = \{y_j,y_k\}$.
    Recall~\eqref{eq:hatpi-def} makes $\omega^{(0)}_u(y_i) = \mu_X(\{y_i\})$,
    and so,
    \[
        \frac{\omega^{(0)}_u(y_i)+\omega^{(0)}_u(y_j) + \omega^{(0)}_u(y_k)}{2} - \omega^{(0)}_u(y_i) \geq \frac{\omega^{(0)}_u(y_i) + \omega^{(0)}_u(y_i)}{2} - \omega^{(0)}_u(y_i) = 0.
    \]
   
    The next claim shows that the transformations encode the relevant projections as cycles or nonempty cuts while preserving their weights and intersections with \(R_u(A)\).
    \begin{claim}\label{claim:phi}
        Suppose $(E_u,M_u)$ is graphic or cographic,
        and let $Z\in\mathscr Z(M_u)$ be a nonempty cycle satisfying $Z\cap F_u =\varnothing$ and $|Z\cap S_{u,X}| \leq 1$ for all $X\in\mathcal X_3(u)$.
        Then, $\Phi_u(Z)$ is a nonempty even subgraph of $G_u^{(A)}$ if $M_u$ is graphic, or a nonempty cut if cographic, such that $\omega^{(A)}_u(\Phi_u(Z)) = \omega^{(0)}_u(Z)$ and $\Phi_u(Z)\cap R_u(A)=Z\cap R_u(A)$.
        Furthermore, if $Z\in\mathscr C(M_u)$ and $M_u$ is graphic then $\Phi_u(Z)$ is also a cycle of $G^{(A)}_u$.
    \end{claim}
    \begin{croof}
        Fix nonempty $Z\in\mathscr Z(M_u)$, and first note that 
        because every \(X\in\mathcal X_3(u)\) has \(u\) as its unique selected neighbor, the separator $S_{u,X}$ is disjoint from $R_u(A)$ as $R_u(A)$ only records separator choices on paths toward other selected vertices.
        Therefore, since the exchanges defining~\eqref{eq:Phi-def} only occur within such separators, $\Phi_u(Z)\cap R_u(A)=Z\cap R_u(A)$ follows, as well as $\Phi_u(Z)\cap F_u = \varnothing$ since $Z\cap F_u =\varnothing$.

        Now, because $|Z\cap S_{u,X}|\leq 1$ for all $X\in\mathcal X_3(u)$, any $\Delta$-$Y$ exchanges in $G^{(0)}_u$ meeting $Z$ do not interact.
        Each exchange preserves a nonempty even subgraph in the graphic case and a nonempty cut in the cographic case. If \(Z\) is a circuit and \(\mathcal M_u\) is graphic, $\Phi_u(Z)$ is a cycle.
        Applying $\Phi_u(Z)\cap F_u = \varnothing$ then shows $\Phi_u(Z)$ persists into $G^{(A)}_u$ with these properties.

        Continuing on, observe that for any $X\in\mathcal X_3(u)$ and separator $S_{u,X}$ meeting $Z$, our assumptions make $|Z\cap S_{u,X}| = 1$.
        So, let $S_{u,X} = \{y, z_1, z_2\}$ and $Z\cap S_{u,X} = \{y\}$.
        Then,
        \begin{multline*}
            \omega^{(A)}_u(S_{u,X}\setminus\{y\}) = \omega^{(A)}_u(z_1) + \omega_u^{(A)}(z_2)\\ = \left(\frac{\omega_u^{(0)}(S_{u,X})}{2}-\omega_u^{(0)}(z_1)\right)
+\left(\frac{\omega_u^{(0)}(S_{u,X})}{2}-\omega_u^{(0)}(z_2)\right) =\omega^{(0)}_u(y).
        \end{multline*}
        Using $\Phi_u(Z)\cap F_u = \varnothing$, this implies
        \begin{equation*}
            \begin{split}
                \omega^{(A)}_u(\Phi_u(Z)) &= \omega^{(0)}_u\left(\Phi_u(Z)\setminus \bigcup_{X\in\mathcal X_3(u)} S_{u,X}\right) + \sum_{(y,S_{u,X}):y\in Z\cap S_{u,X},X\in\mathcal X_3(u)} \omega^{(A)}_u(S_{u,X}\setminus\{y\})\\
                &= \omega^{(0)}_u\left(Z\setminus \bigcup_{X\in\mathcal X_3(u)} S_{u,X}\right) + \sum_{(y,S_{u,X}):y\in Z\cap S_{u,X},X\in\mathcal X_3(u)} \omega^{(0)}_u(y)\\
                &=\omega^{(0)}_u(Z).\qedhere
            \end{split}
        \end{equation*}
    \end{croof}
    Now we may define $C\mapsto (D_u^{(A)}(C))_{u\in \mathcal S(A)}$, with $C\in\sigma^{-1}(A)$, via $D_u^{(A)}(C)\triangleq \Phi_u(C_u)$.
    This is injective:
    Fix two circuits $C^{(1)},C^{(2)}\in\sigma^{-1}(A)$ such that $D^{(A)}_u(C^{(1)})=D^{(A)}_u(C^{(2)})$ for all $u\in\mathcal S(A)$.
    Rearranging~\eqref{eq:Phi-def} shows
    \begin{equation}\label{eq:cu-reverse}
        C_u =D^{(A)}_u(C)\symdiff\left(\bigsymdiff_{X\in\mathcal X_3(u):C_u\cap S_{u,X}\neq\varnothing} S_{u,X}\right),\quad\forall C\in\sigma^{-1}(A).
    \end{equation}
    Furthermore,~\eqref{eq:Phi-def} and Proposition~\ref{prop:projections}.2 also show
    \[
        D_u^{(A)}(C)\cap S_{u,X}=\begin{cases}
            S_{u,X}\setminus\{y\},\quad& \text{if }(\exists y\in E_u)\;C_u\cap S_{u,X}=\{y\},\\
            \varnothing,\quad& \text{otherwise},\\
        \end{cases}
    \]
    and so $C_u\cap S_{u,X}\neq\varnothing$ if and only if $D_u^{(A)}(C)\cap S_{u,X}\neq\varnothing$.
    Applying this observation to~\eqref{eq:cu-reverse} shows $C^{(1)}_u=C^{(2)}_u$ for all $u\in\mathcal S(A)$.
    This implies
    \[
        \partial_X C^{(1)}_X = \bigcup_{u\in N(X)} C^{(1)}_u\cap S_{u,X} = \bigcup_{u\in N(X)} C^{(2)}_u\cap S_{u,X} = \partial_X C^{(2)}_X,
    \]
    so Lemma~\ref{lem:C-assignment}.1 and Lemma~\ref{lem:bd=>equal} show $C^{(1)}_X = C^{(2)}_X$ for all $X\in\mathcal X$ as well.
    Hence,
    \[
        \left(\bigsymdiff_{u\in\mathcal S(A)} C^{(1)}_u\right)\symdiff\left(\bigsymdiff_{X\in\mathcal X}C^{(1)}_X\right) = \left(\bigsymdiff_{u\in\mathcal S(A)} C^{(2)}_u\right)\symdiff\left(\bigsymdiff_{X\in\mathcal X}C^{(2)}_X\right),
    \]
    and, as a result, the factorization of a circuit as supplied by Proposition~\ref{prop:projections} shows $C^{(1)}=C^{(2)}$.

    Select some $C\in\sigma^{-1}(A)$.
    Observe $M_u\cong \mathcal R_{10}$ implies $C_u = D_u^{(A)}(C)$ and $\omega^{(A)}_u(C_u) = \omega^{(0)}_u(C_u)$
    since $\mathcal R_{10}$ matroids only participate in 1- and 2-sums in $T$ by Proposition~\ref{prop:foundation}.3.
    Moreover, by Claim~\ref{claim:phi} we have that $D_u^{(A)}(C)$ is a cycle of $G^{(A)}_u$ when $M_u$ is graphic, and a nonempty cut when cographic.
    Hence 2. follows.
    Continuing on, since $R_u(C) = R_u(A)$ and $R_u(C)\subseteq C_u$, we apply $\Phi_u(C_u)\cap R_u(A)=C_u\cap R_u(A)$ as in Claim~\ref{claim:phi} to see $R_u(C)\subseteq D^{(A)}_u(C)$.
    Claims~\ref{claim:phi} and~\ref{claim:hatpi} together show
    \[
        \omega_u^{(A)}(D_u^{(A)}(C)) = \omega_u^{(A)}(\Phi_u(C_u)) = \omega^{(0)}_u(C_u)=w(C_{B_u(C)}\cap E),
    \]
    and so we've shown 1.
    What remains is 3.

    Let $D$ be a cycle or nonempty cut of $G_u^{(A)}$, corresponding to the cases where $M_u$ is graphic or cographic, satisfying $D\cap R_u(A)=\varnothing$.
    Since $D\cap F_u=\varnothing$, we may regard $D$ as a cycle or nonempty cut of $G_u^{(1)}$, respectively: deleting $F_u$ preserves cycles avoiding $F_u$, while cuts of $G_u^{(1)}/F_u$ correspond to cuts of $G_u^{(1)}$ avoiding $F_u$.
    For each $X\in\mathcal X_3(u)$, the separator $S_{u,X}$ forms a degree-three star in $G_u^{(1)}$ in the graphic case and a triangle in the cographic case.
    Consequently, $|D\cap S_{u,X}|\in\{0,2\}$. 
    Define
    \[
        Z\triangleq D\symdiff\left(\bigsymdiff_{X\in\mathcal X_3(u):D\cap S_{u,X}\neq\varnothing} S_{u,X}\right),
    \]
    i.e. $Z =\Phi_u(D)$.
    This reverses each exchange by replacing a pair of separator elements with the remaining element.
    The resulting set is a cycle of $G_u^{(0)}$ in the graphic case and a cut of $G_u^{(0)}$ in the cographic case.
    Moreover, $Z\ne\varnothing$, since $D\ne\varnothing$ and each replacement preserves nonemptiness. Thus $Z\in\mathscr Z(M_u)\setminus\{\varnothing\}$, with $D=\Phi_u(Z)$ and $|Z\cap S_{u,X}|\le1$ for every $X\in\mathcal X_3(u)$.
    The exchanges affect only separators disjoint from $R_u(A)\cup F_u$, so $Z\cap\bigl(R_u(A)\cup F_u\bigr)=\varnothing$.
    Importantly, this implies $Z$ avoids every separator $S_{u,v}$ for $v\in N^*_u(\mathcal S(A))$. 

    Let $\mathcal X_1(u)$ denote the components of $T-\mathcal S(A)$ incident only to $u$.
    For every $X\in\mathcal X_1(u)$, we have $|Z\cap S_{u,X}|\le1$.
    This follows from the reversed exchanges for three-element separators and is immediate for separators of size zero or one.
    If this intersection is a singleton, Claim~\ref{claim:y-is-boundary} guarantees a cycle realizing that boundary; if it is empty, the empty cycle does so as well.
    Therefore, we may choose a minimizing cycle $Y^{(X)}\in\mathscr Z(M_X)$ satisfying
    \begin{align}\label{eq:ZX}
        \partial_XY^{(X)}=Z\cap S_{u,X}
        &&\text{and}&&
        w(Y^{(X)}\cap E)=\mu_X(Z\cap S_{u,X}).
    \end{align}
    Note $\mu_X(\varnothing)=0$. Therefore, whether $Z\cap S_{u,X}$ is empty or a singleton,~\eqref{eq:hatpi-def} gives $\mu_X(Z\cap S_{u,X})=\omega_u^{(0)}(Z\cap S_{u,X})$.
    Moreover, since $Z$ avoids all other separators, these pairwise disjoint intersections partition $Z\setminus E$, i.e. $\bigcup_{X\in\mathcal X_1(u)}Z\cap S_{u,X} = Z\setminus E$, and
    the symmetric difference
    \[
        Z'\triangleq Z\symdiff\bigsymdiff_{X\in\mathcal X_1(u)} Y^{(X)}
    \]
    is a cycle of $\mathcal M$.
    Observe that the components of $\mathcal X_1(u)$ induce disjoint ground sets in their matroids by Definition~\ref{def:udt}.2, so the intersections $(Y^{(X)}\cap E)_{X\in\mathcal X_1(u)}$ are mutually disjoint.
    Since separator elements cancel in the difference defining $Z'$, it follows that
    \begin{equation}\label{eq:dju}
        Z' = (Z\cap E)\uplus\left(\biguplus_{X\in\mathcal X_1(u)} Y^{(X)}\cap E\right),
    \end{equation}
    and so by~\eqref{eq:ZX} and our prior observations,
    \begin{equation}\label{eq:wZ-to-piad}
        \begin{split}
            w(Z') &= w(Z\cap E) + \sum_{X\in\mathcal X_1(u)}w(Y^{(X)}\cap E)\\
                 &= w(Z\cap E) + \sum_{X\in\mathcal X_1(u)}\mu_X(Z\cap S_{u,X})\\
                 &= \omega^{(0)}_u(Z\cap E) + \sum_{X\in\mathcal X_1(u)}\omega_u^{(0)}(Z\cap S_{u,X})\\
                 &= \omega^{(0)}_u(Z).
        \end{split}
    \end{equation}
    Recalling $Z=\Phi_u(D)$, it follows by Claim~\ref{claim:phi} that $w(Z') = \omega^{(A)}_u(D)$.
%

    To see that $Z'\neq\varnothing$, recall that $Z\neq\varnothing$.
    If $Z\cap E\neq\varnothing$, then~\eqref{eq:dju} gives the result.
    Otherwise, choose $y\in Z$.
    Since $Z$ avoids every separator leading toward another selected vertex, $y$ lies in $S_{u,X}$ for some $X\in\mathcal X_1(u)$.
    Thus $Z\cap S_{u,X}=Y^{(X)}\cap S_{u,X}=\{y\}$.
    In a $2$-sum the shared element is a nonloop, and in a $3$-sum it belongs to a three-element circuit. 
    Thus $Y^{(X)}$ contains an element other than $y$. 
    Since $u$ is the only neighbor of $X$, we have $E_X\setminus E=S_{u,X}$, so this other element belongs to $Y^{(X)}\cap E$ and hence to $Z'$. Therefore $Z'$ is a nonempty cycle of $\mathcal M$, giving $w(Z')\geq g_w$.
    By $w(Z')=\omega^{(A)}_u(D)$ we conclude 3.
\end{proof}

\subsection{Combining Local Counting Estimates}\label{subsec:combining}

In this section we finish the proof of the main result.
Our argument distinguishes circuits according to the number of selected vertices.
When there is only one selected vertex, \(R_u(C)=\varnothing\).
Ao Lemma~\ref{lem:mapping}.3 ensures that every cycle or nonempty cut in the corresponding graph has weight at least $g_w$.
This allows for a direct application of classical bounds of Karger and Subramanian.
With multiple selected vertices, however, the local graphs may contain arbitrarily lighter cycles and cuts and Lemma~\ref{lem:mapping}.3 only guarantees the local cycles and cuts avoiding $R_u(C)$ have weight at least $g_w$.
We address this in Lemma~\ref{lem:heavy-D-bound} by establishing a counting estimate under precisely this assumption.

We begin by recalling the classic bounds.
We note that the bound for cuts as we describe it is shown exactly in \cite{KS96}.
Moreover, \cite{subramanian1995polynomial} as published only concerns half-integral $\alpha$.
Applying this directly would incur extra polynomial factors in our bound, so we provide an argument extending Subramanian's result to all real $\alpha\geq 1$ in Appendix~\ref{app:sub}.
\begin{lemma}[\cite{karger1993global,KS96,subramanian1995polynomial}]\label{lem:other-bound}
    Let $G$ be a graph on at most $m\geq 1$ edges, $\omega:\mathcal E(G)\to\mathbb{R}_{\geq 0}$ a nonnegative weight function, $g_w>0$,
    and $\mathscr D$ be either the family of cycles or the family of nonempty cuts.
    For any choice of $\mathscr{D}$, if $\omega(D) \geq g_w$ for every $D \in \mathscr D$ then for every real $\alpha\geq 1$,
    \[
        \big|\big\{D\in\mathscr D\mid \omega(D) \leq \alpha g_w\big\}\big|\leq (2m)^{2\alpha}.
    \]
\end{lemma}
The next lemma bounds a weighted sum of cycles or cuts, rather than their number below a threshold. In Claim~\ref{claim:1}, products of these local estimates add the weights in the exponent and recover the original circuit weight. The prescribed set \(R\) will be \(R_u(A)\).
The proof uses ideas from \cite{subramanian1995polynomial} and Karger's random contraction method \cite{karger1993global}, building on applications in \cite{gurjar2021isolating} to counting cycles and cuts containing a single fixed edge.
%
\begin{lemma}\label{lem:heavy-D-bound}
    Let $G$ be a graph on at most $m\geq 1$ edges, $\omega:\mathcal E(G)\to\mathbb{R}_{\geq 0}$ a nonnegative weight function, $g_w>0$,
    and $\mathscr D$ be either the family of cycles or the family of nonempty cuts.
    For any choice of $\mathscr{D}$, there exists a universal constant $K\geq 1$ such that, for every nonempty $R\subseteq \mathcal E(G)$, if $\omega(D) \geq g_w$ for all $D\in \mathscr D$ such that $R\cap D =\varnothing$,
    \[
        \sum_{D \in \mathscr D: D\supseteq R} (Km)^{-2\omega(D)/g_w}\leq K^{|R|}(|R|)!.
    \]
\end{lemma}
\begin{proof}
    First select a cycle $D\supseteq R$ and perform the following procedure:
    Mark every edge in $D\cap R$, and traverse $D$ in a fixed order beginning at some fixed edge in $R$.
    While traversing, mark an additional edge whenever the accumulated weight since the preceding mark first reaches $g_w/2$, and reset the accumulated weight after every mark.
    Observe $k\leq 2\omega(D)/g_w$ if $k$ additional edges are marked since each additional mark accounts for a disjoint portion of the cycle of weight at least $g_w/2$.
    Moreover, every path strictly between consecutive marked edges avoids $R$ and has weight strictly less than $g_w/2$.
    This implies $D$ is determined entirely by its marked edges, their ordering in the traversal, and their orientation in the traversal.
    Specifically, if two distinct cycles possess the same set of marked edges with the same orientation and ordering, then there exists a pair of vertices incident to consecutive marked edges with distinct paths between them in the two cycles.
    See that the union of these two paths would contain a cycle avoiding \(R\) of weight less than \(g_w\), which contradicts our assumptions.
    Let $N_k$ be the number of cycles containing $R$ with exactly $k$ additional marked edges after performing this procedure.
    Observe that
    \[
        N_k \leq \binom{m}{k}(|R|+k)!2^{|R|+k}\leq \frac{m^k}{k!}(|R|+k)!2^{|R|+k}\leq m^k(|R|)!\binom{|R|+k}{|R|} 2^{|R|+k}\leq m^k(|R|)!4^{|R| +k},
    \]
    so, if $\mathscr D$ gives the cycles,
    \[
        \sum_{D \in \mathscr D: D\supseteq R} (Km)^{-2\omega(D)/g_w} \leq \sum_{k\geq 0}^\infty N_k (Km)^{-k} \leq 4^{|R|}(|R|)!\sum_{k\geq 0}^{\infty}\left(\frac{4}{K}\right)^{k}.
    \]
    Finally, $4^{|R|}(|R|)!\sum_{k\geq 0}^{\infty}\left(\frac{4}{K}\right)^{k}\leq K^{|R|}(|R|)!$ for all $K\geq 8$, so we have this case.

    Now we examine the case where $\mathscr D$ gives nonempty cuts.
    We may assume that $G$ is connected. 
    Indeed, if $G$ is disconnected, choose one vertex from each connected component and merge these vertices into a single common vertex, retaining all edges and their weights.
    This preserves all cuts and gives a number of vertices $n_0\leq m+1$. 
    Next, assign each positive-weight edge $y\notin R$ an independent exponential clock of rate $\omega(y)$.
    Contract edges when their clocks ring, ignoring loops and leaving zero-weight edges uncontracted.
    At time $\tau \triangleq g_w^{-1}\cdot2\log m$ choose a uniformly random bipartition of the remaining vertices.
    Output its cut if it contains $R$, and otherwise fail.
    If an edge of $R$ becomes a loop, the experiment can immediately fail.
    Fix a target cut $D\supseteq R$, and condition on the event $\mathcal A_D$ that no positive-weight edge of $D\setminus R$ has its clock ring before the stopping time $\tau$.
    Since the clocks are independent, the conditional process is equivalent to contracting only edges outside $D\cup R$ according to their original exponential clocks.
    For each $t\in[0,\tau]$, let $N_t$ denote the number of vertices remaining after the contractions occurring by time $t$.
    At any state, at most $2|R|$ vertices are incident to $R$.
    Since the graph remains connected, every other vertex has a nonempty singleton cut avoiding $R$ of weight at least $g_w$, so total edge weight outside $R$ is at least $(N_t-2|R|)g_w/2$.
    Then conditioning on $\mathcal A_D$ leaves a rate of at least $(N_t-2|R|)g_w/2 - \omega(D)$.
    Since each contraction decreases the vertex count by one, we obtain
    \[
        \frac{d}{dt}\mathbb{E}[N_t\mid \mathcal A_D] \leq -\frac{(\mathbb{E}[N_t\mid \mathcal A_D]-2|R|)g_w}{2} + \omega(D).
    \]
    Multiplying by \(e^{g_wt/2}\) and applying the product rule gives
    \begin{equation*}
        \begin{split}
            0 &\geq e^{g_wt/2}\cdot\frac{d}{dt}\mathbb{E}[N_t\mid \mathcal A_D] + \frac{e^{g_wt/2}g_w}{2}\left(\mathbb{E}[N_t\mid \mathcal A_D]-2|R| - \frac{2\omega(D)}{g_w}\right)\\
              &= \frac{d}{dt}\left[e^{g_wt/2}\cdot\left(\mathbb{E}[N_t\mid \mathcal A_D]-2|R| - \frac{2\omega(D)}{g_w}\right)\right].
        \end{split}
    \end{equation*}
    Integrating over \(t\in[0,\tau]\) and rearranging shows
    \begin{equation*}
        \begin{split}
            \mathbb{E}[N_\tau\mid \mathcal A_D]&\leq 2|R| +\frac{2\omega(D)}{g_w} + \left(n_0 - 2|R| - \frac{2\omega(D)}{g_w}\right)e^{-g_w\tau/2}\\
            &= 2|R| +\frac{2\omega(D)}{g_w} + m^{-1}\cdot\left(n_0 - 2|R| - \frac{2\omega(D)}{g_w}\right)\\
            &\leq 2|R| + \frac{2\omega(D)}{g_w} + 1.
        \end{split}
    \end{equation*}
    Given a surviving graph with $N_\tau$ vertices, a uniformly random bipartition produces the cut $D$ with probability $2^{1-N_\tau}$.
    So, by applying Jensen's inequality,
    \begin{equation}\label{eq:jensens}
        \Pr[\text{output }D\mid\mathcal A_D] = \mathbb{E}\left[2^{1-N_\tau}\mid\mathcal A_D\right] \geq 2^{-2|R| -2\omega(D)/g_w}.
    \end{equation}
    By independence,
    \[
        \Pr[\mathcal A_D]=e^{-\tau\omega(D\setminus R)} =m^{-2\omega(D\setminus R)/g_w}\geq m^{-2\omega(D)/g_w},
    \]
    and so combining with~\eqref{eq:jensens} implies
    \begin{equation*}
            \Pr[\text{output }D] \geq 2^{-2|R| -2\omega(D)/g_w}\cdot m^{-2\omega(D)/g_w}=4^{-|R|}(2m)^{-2\omega(D)/g_w}.
    \end{equation*}
    But then, as the experiment outputs at most one cut, summing over the marginals shows
    \[
        1 \geq \sum_{D\in\mathscr D:D\supseteq R}\Pr[\text{output }D] \geq 4^{-|R|}\sum_{D\in\mathscr D:D\supseteq R}(2m)^{-2\omega(D)/g_w}. 
    \]
    Rearranging, we obtain
    \[
        \sum_{D\in\mathscr D:D\supseteq R}(Km)^{-2\omega(D)/g_w}\leq \sum_{D\in\mathscr D:D\supseteq R}(2m)^{-2\omega(D)/g_w}\leq 4^{|R|}\leq K^{|R|},
    \]
    whenever $K\geq 4$.
\end{proof}

\begin{proof}[Proof of Main Result]
    To simplify the notation let $Q \triangleq \sum_{u\in\mathcal V(T)}|E_u|$ and $K\geq 32$ be a constant large enough to satisfy both Lemmas~\ref{lem:heavy-D-bound} and~\ref{lem:other-bound}.
    Before proceeding, examine a vertex $u\in\mathcal V(T)$ such that $(E_u,\mathcal M_u)\cong \mathcal R_{10}$.
    Note that $K\geq 32$ implies
    \begin{equation}\label{eq:r10-bound}
        \sum_{D\subseteq E_u:D\supseteq R_u(A)}(K|E_u|)^{-2\omega_u^{(A)}(D)/g_w}\leq 2^{10}\leq K^2 \leq K^{|R_u(A)|+1}|R_u(A)|!,
    \end{equation}
    for all signatures $A$ with $u\in\mathcal S(A)$ and $|\mathcal S(A)|\geq 2$ as this ensures $|R_u(A)|\geq 1$.
    In what follows, define $\mathscr D_u^{(A)}$ to be the cycles or nonempty cuts, when $M_u$ is graphic or cographic respectively, of the graph $G^{(A)}_u$ as given by Lemma~\ref{lem:mapping}.
    When $M_u\cong\mathcal R_{10}$ then we simply let $\mathscr D_u^{(A)}$ equal $2^{E_u}$.
    We separately bound the number of $\alpha$-minimum weight circuits with exactly one selected vertex and with at least two selected vertices.
    \begin{claim}\label{claim:1}
    \[\sum_{C\in\mathscr C(\mathcal M):|\mathcal S(C)| \geq 2}\left(18K^4Q\right)^{-2w(C)/g_w}\leq \frac 1 2.\]
    \end{claim}
    \begin{croof}
        Fix a realized signature \(A\), and write \(S=\mathcal S(A)\) and \(t=|S|\ge2\).
        Examine the smallest connected subtree $H$ of $T$ containing $S$.
        Minimality implies the neighbors of $u$ in $H$ are precisely $N_u^*(S)$ for all $u\in S$.
        Let $d_v$ give the degree of any $v\in\mathcal \mathcal V(H)$ in $H$, and observe $d_u=|R_u(A)|$ by~\eqref{eq:proj-sep-single} for all $u\in S$.
        By $|S|\geq 2$ we find the degrees are strictly positive.
        Moreover, all leaves lie in $S$, so $d_v\geq2$ for all $v\in\mathcal V(H)\setminus S$.
        Therefore, by the handshaking lemma,
        \begin{equation}\label{eq:sum-ru-bound}
            \sum_{u\in S} d_u = 2|\mathcal E(H)| - \sum_{v \in \mathcal V(H)\setminus S}d_v
                                            \leq 2(t +|\mathcal V(H)\setminus S| - 1) - 2|\mathcal V(H)\setminus S|
                                = 2t - 2,
        \end{equation}
        and so repeated applications of $a!b!\leq (a+b-1)!$ show
        \begin{equation}\label{eq:prod-ru-bound}
        \begin{split}
            \prod_{u\in S}d_u! \leq \left(1 - t + \sum_{u\in S}d_u\right)! \leq (t-1)!.
        \end{split}
        \end{equation}
        Lemmas~\ref{lem:mapping}.3 and~\ref{lem:heavy-D-bound}, together with~\eqref{eq:r10-bound}, give
        \begin{equation}\label{eq:fkko}
            \sum_{D\in\mathscr D^{(A)}_u:D\supseteq R_u(A)}(K|E_u|)^{-2\omega_u^{(A)}(D)/g_w}\leq K^{d_u+1}d_u!,
        \end{equation}
        and Lemmas~\ref{lem:C-assignment} and~\ref{lem:mapping} also give
        \begin{equation}\label{eq:dfsdf}
            \sum_{u\in \mathcal S(C)}\omega_u^{(A)}(D_u^{(A)}(C))=w(C)\text{ and }\min_{u\in S}\frac{2\omega_u^{(A)}(D_u^{(A)}(C))}{g_w}\geq 1,\quad\forall C\in\sigma^{-1}(A).
        \end{equation}
        Let $L\geq K$ be a constant to be determined later.
        Observe that
        \begin{equation}\label{eq:ll}
            (LQ)^{-\beta} = (K|E_u|)^{-\beta}\left(\frac {K|E_u|}{LQ}\right)^\beta \leq \left(\frac {K|E_u|}{LQ}\right)(K|E_u|)^{-\beta},
        \end{equation}
        for all $\beta\geq1$.
        So~\eqref{eq:dfsdf} and~\eqref{eq:ll} imply
        \begin{equation}\label{eq:run0}
            \begin{split}
                \sum_{C \in \sigma^{-1}(A)}\left(L Q\right)^{-2w(C)/g_w}&= \sum_{C \in \sigma^{-1}(A)}\prod_{u\in S}(LQ)^{-2\omega_u^{(A)}(D_u^{(A)}(C))/g_w}\\
                                                                        &\leq \sum_{C \in \sigma^{-1}(A)}\prod_{u\in S}\frac{K|E_u|}{LQ}(K|E_u|)^{-2\omega_u^{(A)}(D_u^{(A)}(C))/g_w}\\
                                                                        &= \left(\prod_{u\in S}\frac{K|E_u|}{LQ}\right)\cdot\sum_{C \in \sigma^{-1}(A)}\prod_{u\in S}(K|E_u|)^{-2\omega_u^{(A)}(D_u^{(A)}(C))/g_w}.
            \end{split}
        \end{equation}
        By Lemma~\ref{lem:mapping}, \(C\mapsto(D_u^{(A)}(C))_{u\in S}\) is injective, with image contained in the Cartesian product of the families \(\{D\in\mathcal D_u^{(A)}:D\supseteq R_u(A)\}\) over $u\in S$. Since all summands are nonnegative, enlarging the image to this Cartesian product, applying distributivity, and applying~\eqref{eq:fkko} gives
        \[
            \sum_{C \in \sigma^{-1}(A)}\prod_{u\in S}(K|E_u|)^{-2\omega_u^{(A)}(D_u^{(A)}(C))/g_w}
            \leq \prod_{u\in S} \sum_{D\in\mathscr D^{(A)}_u:D\supseteq R_u{(A)}} (K|E_u|)^{-2\omega_u^{(A)}(D)/g_w}\leq \prod_{u\in S} K^{d_u+1}d_u!.
        \]
        And so, via~\eqref{eq:sum-ru-bound} and~\eqref{eq:prod-ru-bound},
        \begin{equation}\label{eq:run1}
            \begin{split}
                \left(\prod_{u\in S}\frac{K|E_u|}{LQ}\right)\cdot\sum_{C \in \sigma^{-1}(A)}\prod_{u\in S}(K|E_u|)^{-2\omega_u^{(A)}(D_u^{(A)}(C))/g_w}&\leq \prod_{u\in S}\left(\frac {K|E_u|}{LQ}\cdot  K^{d_u+1}d_u!\right)\\
                                                                         &\leq \frac{K^{4t-2}(t-1)!}{L^t}\cdot\prod_{u\in S}\frac{|E_u|}{Q}.
            \end{split}
        \end{equation}
            Now, by~\eqref{eq:proj-sep-single}, there are at most $3^{d_u}$ choices of $R_u(A)$, and so by~\eqref{eq:sum-ru-bound} there are at most $3^{2t-2}$ possible signatures which could have a set of selected vertices equal to $S$.
        Hence from~\eqref{eq:run0} and~\eqref{eq:run1},
        \begin{equation}\label{eq:asdfi}
            \begin{split}
                \sum_{C\in\mathscr C(\mathcal M):S = \mathcal S(C)}\left(L Q\right)^{-2w(C)/g_w}&=\sum_{A\in\sigma(\mathscr C(\mathcal M)):S=\mathcal S(A)}\sum_{C \in \sigma^{-1}(A)}\left(L Q\right)^{-2w(C)/g_w}\\&\leq\left(\frac{9K^{4}}{L}\right)^t(t-1)!\cdot\prod_{u\in S}\frac{|E_u|}{Q}.
            \end{split}
        \end{equation}
        Finally,
        \[
            t!\sum_{S\subseteq V(T):|S|=t}
            \prod_{u\in S}\frac{|E_u|}{Q}
            \leq
            \left(\sum_{u\in\mathcal V(T)} \frac{|E_u|}{Q}\right)^t
            =1.
        \]
        Dividing the preceding inequality by \(t\), summing~\eqref{eq:asdfi} over all sets \(S\) of size \(t\), and then summing over \(t\ge2\) gives
        \[
            \sum_{C\in\mathscr C(\mathcal M):|\mathcal S(C)|\geq 2}\left(L Q\right)^{-2w(C)/g_w}\leq \sum_{t=2}^\infty \frac 1 t\left(\frac{9K^{4}}{L}\right)^t.
        \]
        Taking $L=18K^4$ then gives the claim.
    \end{croof}
    \begin{claim}\label{claim:2}
        For every $\alpha \geq 1$, $\big|\big\{C\in\mathscr C\bigm\vert |\mathcal S(C)| = 1\text{ and } w(C) \leq \alpha g_w\big\}\big| \leq \left(KQ\right)^{2\alpha}$.
    \end{claim}
    \begin{croof}
        If \(|\mathcal S(C)|=1\), then \(\sigma(C)=\{(u,\varnothing)\}\) for some \(u\in\mathcal V(T)\). Fix \(u\) such that \(\{(u,\varnothing)\}\) is a realized signature; vertices with empty signature classes contribute zero. Lemma 3.3 gives an injection \(C\mapsto D_u(C)\) on this signature class, and Lemmas~\ref{lem:C-assignment}.3 and~\ref{lem:mapping}.1 give \(\omega_u(D_u(C))=w(C)\). Thus,
%
        \[
            C\in\sigma^{-1}(\{(u,\varnothing)\})\text{ and }w(C) \leq \alpha g_w\implies \omega_u(D_u(C))\leq \alpha g_w,
        \]
        and so by injectivity of $C\mapsto (D_u(C))_{u\in\{u\}}$,
        \begin{equation}\label{eq:c-in-sig-lt-d}
            |\{C\in\sigma^{-1}(\{(u,\varnothing)\})\mid w(C) \leq \alpha g_w\}| \leq |\{D\in\mathscr D_u^{(\{(u,\varnothing)\})}\mid \omega_u(D)\leq \alpha g_w\}|.
        \end{equation}
        Suppose $(E_u,\mathcal M_u)$ is graphic or cographic.
        Since $C\in\sigma^{-1}(\{(u,\varnothing)\})$ implies $R_u(C)=\varnothing$, Lemma~\ref{lem:mapping}.3 makes $\omega_u(D)\geq g_w$ for all $D \in \mathscr D_u^{(\{(u,\varnothing)\})}$.
        Therefore, we may apply Lemma~\ref{lem:other-bound} to obtain
        \begin{equation}\label{eq:applying-other-bound}
            |\{D\in\mathscr D_u^{(\{(u,\varnothing)\})}\mid \omega_u(D)\leq \alpha g_w\}| \leq (K|E_u|)^{2\alpha},
        \end{equation}
        when $M_u$ is graphic or cographic.
        If instead $M_u\cong\mathcal R_{10}$, then $|E_u| = 10$ and $\mathscr D_u^{(\{(u,\varnothing)\})} = 2^{E_u}$ has cardinality 1024.
        Since \(|E_u|=10\), \(K\ge32\), and \(\alpha\ge1\), we have \(1024\le(K|E_u|)^{2\alpha}\), proving~\eqref{eq:applying-other-bound} in all cases.
        Combining~\eqref{eq:c-in-sig-lt-d} and~\eqref{eq:applying-other-bound} gives
        \begin{equation*} 
        \begin{split}
            \big|\big\{C\in\mathscr C\bigm\vert |\mathcal S(C)| = 1\text{ and } w(C) \leq \alpha g_w\big\}\big| &= \sum_{u\in\mathcal V(T)} |\{C\in\sigma^{-1}(\{(u,\varnothing)\})\mid w(C) \leq \alpha g_w\}|\\
                                                                                                                                                                                                                      &\leq \sum_{u\in\mathcal V(T)} (K|E_u|)^{2\alpha}.
        \end{split}
        \end{equation*}
        So, because $\sum_{u\in\mathcal V(T)}|E_u|^{2\alpha}\leq Q^{2\alpha}$ we conclude the claim.
    \end{croof}
    We use the prior two claims to finish the proof of the main result.
    Observe $w(C) \leq \alpha g_w$ implies $2w(C)/g_w \leq 2\alpha$.
    So, since $18K^4Q\geq 1$, we find $\left(18K^4Q\right)^{-2w(C)/g_w} \geq \left(18K^4Q\right)^{-2\alpha}$
    for $\alpha$-minimum weight circuits $C$.
    Applying Claim~\ref{claim:1} and rearranging gives
    \[
        \big|\big\{C\in\mathscr C\mid |\mathcal S(C)| \geq 2\text{ and } w(C) \leq \alpha g_w\big\}\big| \leq \frac 1 2\left(18K^4Q\right)^{2\alpha}.
    \]
    Thus, with Claim~\ref{claim:2},
    \begin{equation*}
        \bigl|\big\{C\in\mathscr C(\mathcal M) \mid w(C)\le\alpha g_w\big\}\bigr|\leq  \left(KQ\right)^{2\alpha} + \frac 1 2\left(18K^4Q\right)^{2\alpha}\leq \left(18K^4Q\right)^{2\alpha}.
    \end{equation*}
    Since \(Q\le7m-6\) and \((\mathcal M,w)\) was arbitrary, it follows that $f_\reg(m,\alpha) \leq (126K^4m)^{2\alpha}$.
\end{proof}

\begin{remark}\label{remark:ext}
    The main result extends to binary matroids obtained by 1- and 2-sums of regular matroids and matroids of bounded size. 
    Indeed, refine the regular summands as in Proposition 2.2 and retain the bounded summands as additional nodes; their attachment through only 1- and 2-sums preserves both conclusions of that proposition.
    The arguments of Lemmas~\ref{lem:bd=>equal} and~\ref{lem:C-assignment} extend without change as they only depend on the structure of the decomposition tree.
    Lemma~\ref{lem:mapping} extends by treating any matroids on a bounded-size ground set in the same way we treated $\mathcal R_{10}$ matroids.
    Specifically, take $D_u^{(A)}(C)=C_u$ and $\omega_u^{(A)}=\omega_u^{(0)}$, as in~\eqref{eq:hatpi-def}, for vertices $u\in\mathcal S(A)$ with such matroids $M_u$.
    Then~\eqref{eq:r10-bound} is recovered by increasing the constant $K$ so that $K^2\geq 2^b$ when the largest ground set of one of the bounded summands is sized $b$.
    Remaining pieces of the argument go through unchanged.
\end{remark}

\section{Lower Bounds Conditioned on the Erdős Girth Conjecture}\label{sec:lb}

Now we prove our conditional lower bounds apply to all fixed real $\alpha\geq 1$, under the Erdős Girth Conjecture predicting the existence of graphs whose girth remains large relative to their edge density.

\begin{conjecture}[Erdős Girth Conjecture]
    For every integer $r\ge 1$, there is a constant $K>0$ such that for every sufficiently large $n$ there is a graph $G$ on $n$ vertices satisfying $|E(G)|\ge K n^{1+1/r}$ and girth strictly greater than $2r.$
\end{conjecture}
This is combined with Simonovits' even-cycle supersaturation theorem (first stated in \cite{erdos1984cube}; see \cite{jiang2020supersaturation} for an accessible proof)
showing that sufficiently high edge density guarantees not merely the existence of cycles of a prescribed even length, but a lower bound on their number.
\begin{theorem}[Simonovits' Even-Cycle Supersaturation Theorem \cite{erdos1984cube,jiang2020supersaturation}]
    Let $G$ be a simple graph on $n$ vertices.
    For every fixed integer $q\geq 2$, there exist positive constants $K_1,K_2$ such that if $|\mathcal E(G)|\geq K_1n^{1+1/q}$ then $G$ contains at least $K_2(|\mathcal E(G)|/n)^{2q}$ distinct cycles of length $2q$.
\end{theorem}
Applying supersaturation to the high-girth graphs supplied by the conjecture yields many cycles of length at most \(\alpha\) times the girth, giving lower bounds with exponents arbitrarily close to \(2\alpha\).
\begin{theorem}\label{thm:lower-bound}
    If the Erdős Girth Conjecture holds true, for every real $\alpha\ge1$, $\varepsilon>0$, and sufficiently large $m$ we have $f_\reg(m,\alpha)\ge m^{2\alpha-\varepsilon}$.
\end{theorem}
\begin{proof}[Proof of Theorem~\ref{thm:lower-bound}]
    Fix $\alpha\ge1$ and $\varepsilon>0$. 
    Choose an integer $p\ge3$ with $2/p<\varepsilon$, and set $q=\lfloor\alpha p\rfloor$. 
    Then $q\geq p$ and $2q/p>2\alpha-\varepsilon$.
    Now fix a sufficiently large integer $m$ and take $n=\left\lceil(2m/K)^{(p-1)/p}\right\rceil$.
    The girth conjecture with $r=p-1$ gives an $n$-vertex graph $G$ with at least $Kn^{p/(p-1)}$ edges and girth greater than $2p-2$. 
    Note there is a bipartite subgraph retaining at least half the edges.
    Let $H$ be a subgraph of this bipartite subgraph keeping exactly $m$ edges.
    Every cycle of $H$ is an even cycle of $G$, so the girth of $H$ is $g\geq 2p$. 

    Since $q\ge p$, this implies $m/n^{1+1/q}\to\infty$ as $m\to\infty$.
    Thus, for sufficiently large $m$, Simonovits’s even-cycle supersaturation theorem
    states $H$ contains at least
    \[
        L_1\left(\frac{m}{n}\right)^{2q}\geq L_2m^{2q/p}
    \]
    distinct cycles of length $2q$, where $L_1$ and $L_2$ are constants independent of $m$.
    Every such cycle is a circuit in the graphic matroid represented by $H$ and $2q\leq2\alpha p\leq\alpha g$.
    Consequently $f_\reg(m,\alpha)\geq m^{2\alpha-\varepsilon}(L_2m^{\delta})$ with $\delta\triangleq 2q/p-(2\alpha-\varepsilon)$.
    As $m^{2\alpha-\varepsilon}(L_2m^{\delta})\geq m^{2\alpha-\varepsilon}$ for $m\geq L_2^{-1/\delta}$ we have the result.
\end{proof}

\section{Applications to Lattices and Error-Correcting Codes}\label{sec:app}
Finally we prove Theorems~\ref{thm:codes} and~\ref{thm:lattices}.
Both follow from an elementary lemma that transfers a linear counting exponent from individual objects to sequences constrained by their total weight.

\begin{lemma}\label{lem:ele-cnt}
    Suppose $m\geq2$, $g>0$, $\mathscr D$ is a finite set, and $\omega:\mathscr D\to[g,\infty)$ satisfies
    \begin{equation}\label{eq:cnt-assump}
        \left|\left\{D\in\mathscr D\mid\omega(D)\leq\beta g\right\}\right|\leq(Km)^{2\beta},\quad\forall \beta\geq 1,
    \end{equation}
    for some absolute constant $K$.
    Then, for every fixed real $\alpha\geq1$,
    \[
        \sum_{k=1}^\infty\left|\left\{(D_1,\ldots,D_k)\in\mathscr D^k\Biggm\vert\sum_{i=1}^k\omega(D_i)\leq\alpha g\right\}\right|\leq m^{2\alpha+o(1)}.
    \]
\end{lemma}
\begin{proof}
    Fix $\varepsilon>0$ and let $Z_\varepsilon\triangleq\sum_{D\in\mathscr D}(Km)^{-(2+\varepsilon)\omega(D)/g}$.
    Also define $f(t)\triangleq |\{D\in\mathscr D\mid \omega(D)\leq tg\}|$.
    Then observe that
    \begin{equation*}
        \begin{split}
            (Km)^{-(2+\varepsilon)\omega(D)/g} &= (2+\varepsilon)\ln (Km)\cdot\int_{\omega(D)/g}^\infty e^{-t(2+\varepsilon)\ln (Km)}dt\\
                                                 &= (2+\varepsilon)\ln (Km)\cdot\int_{1}^\infty \mathbf{1}\{\omega(D)\leq tg\}\cdot(Km)^{-t(2+\epsilon)}dt,
        \end{split}
    \end{equation*}
    where the last equality uses $\omega(D)/g\geq 1$.
    By interchanging sum and integration and applying~\eqref{eq:cnt-assump},
    \begin{equation*}
        \begin{split}
            Z_\varepsilon &= (2+\varepsilon)\ln (Km)\cdot\int_{1}^\infty f(t)\cdot(Km)^{-t(2+\varepsilon)}dt\\
                          &\leq (2+\varepsilon)\ln (Km)\cdot\int_{1}^\infty (Km)^{-\varepsilon t}dt\\
                          &=\frac{2+\varepsilon}{\varepsilon}(Km)^{-\varepsilon}.
        \end{split}
    \end{equation*}
    This makes $Z_\varepsilon \to 0$ as $m\to\infty$.
    Next, any $k$-length tuple $(D_1,\ldots,D_k)$ satisfying $\sum_i\omega(D_i)\leq\alpha g$ gives
    \[
        1\leq(Km)^{(2+\varepsilon)\alpha}\prod_{i=1}^k(Km)^{-(2+\varepsilon)\omega(D_i)/g}.
    \]
    Since \(Z_\varepsilon\le1/2\) for sufficiently large \(m\), we have
    \begin{equation*}
    \begin{split}
        \sum_{k=1}^\infty
    \left|
    \left\{
        (D_1,\ldots,D_k)\in\mathscr D^k\Biggm\vert \sum_{i=1}^k\omega(D_i)\leq\alpha g
    \right\}
    \right|
                         &\leq
    (Km)^{(2+\varepsilon)\alpha}
    \sum_{k\geq1}
    \sum_{(D_i)\in\mathscr D^k}
    \prod_{i=1}^k
    (Km)^{-(2+\varepsilon)\omega(D_i)/g}\\
    &=
    (Km)^{(2+\varepsilon)\alpha}
    \sum_{k\geq1}Z_\varepsilon^k\\
    &=
    (Km)^{(2+\varepsilon)\alpha}
    \frac{Z_\varepsilon}{1-Z_\varepsilon}\\
    &\leq (Km)^{(2+\varepsilon)\alpha}.
    \end{split}
    \end{equation*}
    Since $\varepsilon>0$ is arbitrary and $K,\alpha$ are fixed, the exponent can be made arbitrarily close to $2\alpha$ for all sufficiently large $m$.
    Hence the lemma.
\end{proof}

\begin{proof}[Proof of Theorem~\ref{thm:codes}]
    Let $A$ be a parity-check matrix for the code. The support of each nonzero codeword is a disjoint union of circuits of $M[A]$, and the circuit sizes sum to its Hamming weight.
    The distance equals the girth of $M[A]$, and so
    when we take $\mathscr D$ to be the set of circuits, $g=d$, and $\omega(C)=|C|$
    we satisfy the premise of Lemma~\ref{lem:ele-cnt} by our main result.
    Applying Lemma~\ref{lem:ele-cnt} then gives the theorem.
\end{proof}

Now, let $A$ be an $n\times m$ totally unimodular matrix.
A \emph{signed circuit vector} is a nonzero primitive vector (its nonzero entries have greatest common divisor one) $z\in \mathcal L(A)$ having inclusion-minimal support.
Each circuit of the column matroid $M[A]$ is the common support of two opposite signed circuit vectors with entries 1, -1, or 0 (see \cite{white1987combinatorial,onn2010nonlinear}).
For $x,z\in\mathbb Z^m$, write $x\sqsubseteq z$ if
\[
    x_iz_i\geq0\text{ and }|x_i|\leq|z_i|,\quad\forall i=1,\ldots,m.
\]
When $x\sqsubseteq z$, we say $x$ is \emph{conformal} to $z$.
It is known that nonzero vectors in the nullspace lattice may be written as a conformal sum of signed circuit vectors.
\begin{lemma}[\cite{onn2010nonlinear}]\label{lem:conformal-decomp}
    Let $A$ be a totally unimodular matrix.
    Then, for all nonzero $z \in \mathcal L(A)$ there exists an integer $k$ and a sequence $(x^{(1)},\ldots,x^{(k)})$ of vectors such that $z=\sum_{j=1}^kx^{(j)}$, where each $x^{(j)}$ is a signed circuit vector of $A$ and is conformal to $z$.
\end{lemma}
Notice that conformality forces the coordinates $\{x^{(j)}_i\mid j=1,\ldots,k\}$, with fixed $i$, to have the same sign.
And so, there are no cancellations in the sum defined by the previous lemma.

\begin{proof}[Proof of Theorem~\ref{thm:lattices}]
    Assume $\mathcal L(A)\neq\{0\}$, and fix $p\geq1$ and $\alpha\geq 1$.
    Let $g$ be the girth of $M[A]$, $\mathscr D$ be the finite set of signed circuit vectors of $A$,
    and $\omega(x) = |\supp(x)|$ for all $x\in\mathscr D$.
    
    Select nonzero $z \in \mathcal L(A)$ and a conformal decomposition $z=\sum_{j=1}^kx^{(j)}$ into signed circuit vectors as guaranteed by Lemma~\ref{lem:conformal-decomp}.
    Since the entries of the signed circuit vectors are in $\{-1,1,0\}$ and no cancellations occur in the sum,
    \begin{equation}\label{eq:omg-decomp}
        \sum_{j=1}^k \omega(x^{(j)}) = \sum_{j=1}^k\sum_{i=1}^m|x^{(j)}_i| =\sum_{i=1}^m|z_i| \leq \sum_{i=1}^m|z_i|^p = \lVert z\rVert^p_p.
    \end{equation}
    But because $\supp(x^{(j)})\in\mathscr{C}(M[A])$ for each $j=1,\ldots,k$, it follows that each vector in the decomposition satisfies $\omega(x^{(j)})\geq g$.
    Therefore, by~\eqref{eq:omg-decomp} we find $\lVert z\rVert^p_p\geq g$ for all nonzero $z \in \mathcal L(A)$.
    A signed circuit vector with support size $g$
    attains equality, so $\lambda_p^p=g$. The main result
    implies
    \[
        |\{u\in\mathscr D\mid\omega(u)\leq\beta g\}|=2\left|\left\{C\in\mathscr C(M[A])\bigm\vert|C|\leq\beta g\right\}\right|\leq2(Km)^{2\beta},
    \]
    for every $\beta\geq1$,
    and so $\mathscr D$, $\omega$, and $g$ satisfy the premise of Lemma~\ref{lem:ele-cnt} after increasing the absolute constant $K$.
    Finally, if $0<\|z\|_p\leq\alpha\lambda_p$,~\eqref{eq:omg-decomp} shows its conformal decomposition satisfies
    \[
        \sum_{j=1}^k\omega(x^{(j)})\leq\|z\|_p^p\leq\alpha^p\lambda_p^p=\alpha^p g.
    \]
    Therefore, Lemma~\ref{lem:ele-cnt} gives that the number of nonzero $z\in \mathcal L(A)$ with $\|z\|_p\leq\alpha\lambda_p$ is at most $m^{2\alpha^p + o(1)}$.
\end{proof}

\section*{Acknowledgements}
The author thanks Vijay K. Garg for helpful discussions and support through the completion of this work.
The author was made aware of this problem by the recent conjecture in \cite{khanna_et_al-graphic} and began their investigation by identifying slack in the technical steps of \cite{GV}.
Subsequent observations were given to ChatGPT-5.6, which assisted in exploring and developing proof ideas.
ChatGPT-6/GPT-6 Astra assisted with verifying the consistency of the statements of our lemmas/claims with their usage and identifying logical gaps or typographical mistakes in the written arguments, along with basic grammatical edits.
The author assumes full responsibility for the integrity and accuracy of this work.

\printbibliography

\appendix

\section{Proof of Proposition~\ref{prop:foundation}}\label{app:decomp}

\restatefoundation*
\begin{proof}
    First decompose $\mathcal M$ as 1- and 2-sums of 3-connected matroids; this is always possible due to Seymour \cite{seymour1980decomposition}.
    By Theorem 3 and Proposition C.1 of \cite{aprile2022regular}, any of these 3-connected matroids which are not isomorphic to $\mathcal R_{10}$ can be further refined into a UDT of 3-sums over graphic and cographic matroids with special structure.
    Choose all newly introduced separator elements fresh and mutually distinct.
    We note that their proof is such that when a 3-connected matroid $\mathcal N$ is decomposed, each element of its ground set belongs to the ground set of exactly one of the matroids in the decomposition, while each new separator intersects only the ground sets of matroids participating in its sum (so ground sets corresponding to nonadjacent vertices are disjoint).
    Attach each old \(2\)-sum edge to the unique part containing its separator element; attach \(1\)-sum edges arbitrarily.
    This satisfies the first two items of Definition~\ref{def:udt} as well as 1. of the present result due to Theorem 3 and Proposition C.1 of \cite{aprile2022regular}.

    We now show the third item of Definition~\ref{def:udt}.
    Observe that taking a $k$-sum along a separator preserves the circuits and cocircuits contained in either summand outside that separator. 
    Specifically, let $\mathcal N=\mathcal M_1\oplus_k\mathcal M_2$ and write $S=E_1\cap E_2$ for their separator.
    The definition~\eqref{eq:ksum-def} shows that the cycles contained in $E_1\setminus S$ are unchanged. Indeed, any cycle of $\mathcal M_1$ avoiding $S$ is also a cycle of $\mathcal N$. Conversely, if $Z=Z_1\triangle Z_2$ is a cycle of $\mathcal N$ contained in $E_1\setminus S$, then necessarily $Z_2\subseteq S$. By the defining conditions of the sums, $Z_2$ is either empty or equal to $S$, the latter occurring only for a $3$-sum, where $S$ is also a cycle of $\mathcal M_1$. In either case, $Z$ is a cycle of $\mathcal M_1$. Hence the circuits contained in $E_1\setminus S$ are unchanged.
    Furthermore, since $S$ contains no cocircuit of $\mathcal M_2$ (i.e. $E_2\setminus S$ spans $\mathcal M_2$), there is a basis of $\mathcal M_2$ avoiding $S$.
    For each $s\in S$, its fundamental circuit $C_s$ with respect to this basis satisfies $C_s\cap S=\{s\}$. Thus, for any cycle $Z_1\in\mathcal Z(\mathcal M_1)$, the symmetric difference
    $Z_2\triangleq \mathop{\triangle}_{s\in Z_1\cap S} C_s$ is a cycle of $\mathcal M_2$ satisfying $Z_2\cap S=Z_1\cap S$. Consequently, $Z_1\triangle Z_2\in\mathcal Z(\mathcal N)$, and since $Z_2\cap(E_1\setminus S)=\varnothing$, we have
    \[
        (Z_1\triangle Z_2)\cap(E_1\setminus S)
        =Z_1\cap(E_1\setminus S).
    \]
    Conversely, by the definition of a $k$-sum, every cycle $Z\in\mathcal Z(\mathcal N)$ is the symmetric difference of cycles $Z_1\in\mathcal Z(\mathcal M_1)$ and $Z_2\in\mathcal Z(\mathcal M_2)$ agreeing on $S$. Since $Z_2\cap(E_1\setminus S)=\varnothing$, we have $Z\cap(E_1\setminus S)=Z_1\cap(E_1\setminus S)$.
    So,
    \[
        \{Z\cap(E_1\setminus S)\mid Z\in\mathcal Z(\mathcal N)\}
        =
        \{Z\cap(E_1\setminus S)\mid Z\in\mathcal Z(\mathcal M_1)\}.
    \]
    In a binary matroid, cocycles are precisely the sets meeting every cycle in an even number of elements \cite{white1987combinatorial}.
    Therefore, the previous equation shows the cocycles of $\mathcal M_1$ contained in $E_1\setminus S$
    are exactly the cocycles of \(\mathcal N\) contained in \(E_1\setminus S\). Taking inclusion-minimal nonempty cocycles shows that the cocircuits contained there are unchanged.
    The same argument applies to the other summand $\mathcal M_2$.

    We now verify that the separators in our refinement all initially define valid $k$-sums.
    The $3$-sum edges are already valid within the UDTs given by \cite{aprile2022regular}.
    Each old \(2\)-sum element \(s\) belongs to a unique refined piece and avoids its internal separators. If \(s\) were a loop or coloop there, the preceding observation concerning circuits and cocircuits would give the same property in the reconstructed summand, contradicting the original \(2\)-sum conditions.
    Moreover, the matroid at either endpoint of a \(2\)-sum edge has at least three elements: it is either an unchanged \(2\)-sum factor or participates in an internal \(3\)-sum, in which case it has at least seven.

    We now show that the sums remain valid in every order. Suppose we perform a sum $\mathcal M_1\symdiff\mathcal M_2$ along an edge with separator $S=E_1\cap E_2$. Every remaining separator $R$ is disjoint from $S$ by construction. Hence, as the preceding observation shows that the circuits and cocircuits contained in $R$ are unchanged, the conditions regarding circuits and cocircuits required for $2$- or $3$-sum are preserved. The size conditions also persist, since
    \[
        |E(\mathcal M_1\oplus_k\mathcal M_2)|
        =|E_1|+|E_2|-2|S|
        >\max\{|E_1|,|E_2|\},
    \] 
    where the inequality follows from the size conditions of the sums. Repeating this argument shows that every order of summation is valid.
    For any nonempty \(X\subseteq\mathcal V(T)\) inducing a subtree, every summation order gives ground set \(E_X=\mathop{\triangle}_{u\in X}E_u\).
    We show by induction on \(|X|\) that its cycle family is
    \begin{equation}\label{eq:hyp}
        \left\{
        \bigsymdiff_{u\in X} Z_u
        \;\middle|\;
        (\forall u\in X)\;Z_u\in\mathcal Z( M_u)\text{ and }(\forall (u,v)\in\mathcal E(T[X]))\;
        Z_u\cap S_{u,v}=Z_v\cap S_{u,v}
        \right\},
    \end{equation}
    where $T[X]$ is the subtree induced by $X$.
    The claim is immediate for $|X|=1$. 
    For $|X|\geq 2$, choose any edge $(u,v)\in \mathcal E(T[X])$, and let $X_1$ and $X_2$ be the vertex sets of the two components obtained by deleting $(u,v)$ , with \(u\in X_1\) and \(v\in X_2\). By hypothesis, the cycle spaces of $M_{X_1}$ and $M_{X_2}$ when computed in any summation order follow~\eqref{eq:hyp}.
    So, call the resulting matroids $M_{X_1}$ and $M_{X_2}$.
    The definition~\eqref{eq:ksum-def} of sum is equivalent to obtaining $\mathscr Z(M_{X_1}\symdiff M_{X_2})$ via symmetric differences over cycles $(Z^{(1)},Z^{(2)})\in\mathscr Z(M_{X_1})\times\mathscr Z(M_{X_2})$ satisfying $Z^{(1)}\cap S_{u,v} = Z^{(2)}\cap S_{u,v}$.
    The separator $S_{u,v}$ meets only $E_u$ and $E_v$, however, so by inductive hypothesis there are $Z^{(1)}_u\in\mathscr Z(M_u)$ and $Z^{(2)}_v\in\mathscr Z(M_v)$ such that $Z^{(1)}\cap S_{u,v}=Z^{(1)}_u\cap S_{u,v}$ and $Z^{(2)}\cap S_{u,v} = Z^{(2)}_v\cap S_{u,v}$.
    By~\eqref{eq:ksum-def},
    \begin{equation*}
            \mathscr Z(M_{X_1}\symdiff M_{X_2}) =\{Z^{(1)}\!\!\symdiff Z^{(2)}\mid (Z^{(1)}\!,Z^{(2)})\in\mathscr Z(M_{X_1})\times\mathscr Z(M_{X_2})\text{ and } Z^{(1)}\cap S_{u,v} = Z^{(2)}\cap S_{u,v}\}.
    \end{equation*}
    Applying the inductive hypothesis, this equals
    \begin{equation*}
        \left\{
                                                    \left(\bigsymdiff_{y\in X_1} Z_y\right)
                                                    \symdiff
                                                    \left(\bigsymdiff_{y\in X_2} Z_y\right)
                                                    \;\middle|\;\,
                                                    \begin{array}{@{}l@{}}
                                                    (\forall y\in X_1\cup X_2)\;
                                                        Z_y\in\mathcal Z(M_y),\\[3pt]
                                                    (\forall (y,y')\in\mathcal E(T[X_1])\cup\mathcal E(T[X_2]))\;
                                                        Z_y\cap S_{y,y'}=Z_{y'}\cap S_{y,y'},\\[3pt]
                                                    \text{and }\;Z_u\cap S_{u,v}=Z_v\cap S_{u,v}
                                                    \end{array}
                                                    \right\},
    \end{equation*}
    which simplifies to exactly~\eqref{eq:hyp}.
    So, as matroids are determined by their cycles, performing the sums in any order
    over any nonempty \(X\subseteq\mathcal V(T)\) inducing a subtree gives the same matroid $M_X$.
    Finally, performing the internal $3$-sums first and then the original $1$- and $2$-sums recovers $\mathcal M$, so we've shown the third item of Definition~\ref{def:udt}.

    For 2., each original element occurs once and each separator element twice, so
    \[
        m = \left(\sum_{u\in\mathcal V(T)}|E_u|\right) - 2\cdot\sum_{(u,v)\in\mathcal E(T)}|S_{u,v}|.
    \]
    The minimum ground set size conditions on $k$-sums also  imply $|E_u| > 2|S_{u,v}|$ for all $(u,v)\in\mathcal E(T)$.
    And so, if we root the tree at vertex $r\in\mathcal V(T)$ and let $p(u)$ give the resulting parent vertex of $u$,
    \[
        m = |E_r| +\sum_{u\in\mathcal V(T)\setminus\{r\}}\left(|E_u| - 2|S_{u,p(u)}|\right) \geq |\mathcal V(T)|.
    \]
    Rearranging and using \(|S_{u,v}|\le3\) gives
    \[
         \left(\sum_{u\in\mathcal V(T)}|E_u|\right) = m + 2\cdot\sum_{(u,v)\in\mathcal E(T)}|S_{u,v}| \leq m + 6(|\mathcal V(T)| - 1) \leq m + 6(m-1) = 7m-6.\qedhere
    \]
\end{proof}

\section{Extending Subramanian's Bound to Arbitrary Real $\alpha$}\label{app:sub}

\begin{lemma}
    Let $G$ be a graph on at most $m\geq 1$ edges, $\omega:\mathcal E(G)\to\mathbb{R}_{\geq 0}$ a nonnegative weight function, $g_w>0$, and
    $\mathscr D$ be the family of cycles.
    If $\omega(D)\geq g_w$ for every $D\in \mathscr D$, then for every real $\alpha \geq 1$,
    \[
        \big|\big\{D\in\mathscr D\mid\omega(D)\leq \alpha g_w\big\}\big|\leq (2m)^{2\alpha}.
    \]
\end{lemma}
\begin{proof}
    Our starting observation is that between any two vertices there is at most one path of weight less than $g_w/2$. Otherwise, two distinct such paths would contain a cycle of weight less than $g_w$, a contradiction. So, every nonbacktracking walk of weight less than $g_w/2$ is a path.

    Then, for two oriented edges, write $f\succ e$ when $f$ starts where $e$ ends and is not the reversal of $e$.
    For each oriented edge $e$, consider all finite sequences $\gamma=(e_0,e_1,\ldots,e_k)$ where $e_0=e$, $e_{i+1}\succ e_i$, and
    \[
        \sum_{i=0}^{k-1}\omega(e_i)<\frac{g_w}{2}.
    \]
    Define
    \[
        x_e=
        \sum_{\gamma}
        (2m)^{-\frac{2}{g_w}\sum_{i=0}^{k-1}\omega(e_i)}.
    \]
    The sum includes $\gamma=(e)$, contributing $1$ since the weight of the last edge in the sequence is not summed in the exponent. It has at most $2m$ terms: specifying the last oriented edge $e_k$ uniquely determines the preceding path, since that path has weight less than $g_w/2$. Thus $x_e$ is finite and $x_e\ge1$.
    We claim that
    \begin{equation}\label{eq:xe}
        (2m)^{-2\omega(e)/g_w}\cdot
        \sum_{f\succ e}x_f\le x_e.
    \end{equation}
    Indeed, expanding the left-hand side amounts to prepending $e$ to the sequences counted by $x_f$. Those whose resulting prefix has weight less than $g_w/2$ contribute exactly $x_e-1$, since $e$ is not counted. 
    For the remaining sequences, the prefix after $e$ has weight less than $g_w/2$, so there are at most $2m$ possibilities, determined by their last oriented edge. Each contributes at most
    \[
        (2m)^{-\frac{2}{g_w}(g_w/2)}=\frac1{2m},
    \]
    so their total contribution is at most $1$, proving~\eqref{eq:xe}.

    We now define a random walk. Choose an oriented edge uniformly and traverse it. After traversing $e$, continue along $f\succ e$ with probability
    \[
        p(e,f)\triangleq\frac{x_f}{x_e}(2m)^{-2\omega(e)/g_w},
    \]
    and stop with failure with the remaining probability.
    Notice~\eqref{eq:xe} ensures these are valid probabilities.
    Stop the walk when it first returns to its starting vertex, when its accumulated weight first exceeds
    \[
        \left(\alpha-\frac12\right)g_w,
    \]
    or when it repeats another vertex, in which case it fails. If the walk returns to its starting vertex, output the resulting cycle provided its weight is at most $\alpha g_w$. Otherwise, if it stopped by exceeding the weight threshold, append a path of weight less than $g_w/2$ back to the starting vertex, provided one exists and produces a cycle of weight at most $\alpha g_w$. This closing path is unique by our initial observation. In all other cases, output failure.

    Fix a target cycle $D$ with $\omega(D)\le\alpha g_w$. Choose an orientation of $D$ and let $e_0$ be an oriented edge of $D$ minimizing $e\mapsto x_e$ among its oriented edges.
    Examine the event that the walk starts at $e_0$ and follows $D$ until the stopping rule is triggered.
    Conditioned on this event, we will necessarily output $D$. 
    Indeed, either the walk completes $D$, or it stops after traversing a proper prefix of weight greater than $(\alpha-\tfrac12)g_w$. In the latter case, the complementary path has weight strictly less than
    \[
    \alpha g_w-\left(\alpha-\frac12\right)g_w
    =\frac{g_w}{2},
    \]
    and is therefore the unique closing path.
    Let $e_0,e_1,\ldots,e_k$ be the edges traversed before stopping. Since the stopping threshold was not exceeded before the final edge,
    \[
        \sum_{i=0}^{k-1}\omega(e_i)
        \le\left(\alpha-\frac12\right)g_w.
    \] 
    Moreover, $x_{e_k}\ge x_{e_0}$ by our choice of $e_0$. Hence the probability of following this prescribed sequence, conditional on starting at $e_0$, is
    \[
        \prod_{i=0}^{k-1}p(e_i,e_{i+1})
        =
        \frac{x_{e_k}}{x_{e_0}}(2m)^{-\frac{2}{g_w}\sum_{i=0}^{k-1}\omega(e_i)}
        \ge
        (2m)^{-2(\alpha-1/2)}.
    \]
    Since the probability of selecting $e_0$ initially is at least $1/(2m)$, we conclude
    \[
        \Pr[\text{output }D]
        \ge
        \frac1{2m}(2m)^{-2(\alpha-1/2)}
        =(2m)^{-2\alpha}.
    \]
    Finally, the experiment outputs at most one cycle, so
    \[
        1\ge
        \sum_{\substack{D\in\mathcal D\\\omega(D)\le\alpha g_w}}
        \Pr[\text{output }D]
        \ge
        \left|\{D\in\mathscr D\mid\omega(D)\le\alpha g_w\}\right|
        \cdot(2m)^{-2\alpha}.
    \]
    Rearranging proves the claim.
\end{proof}

\end{document}